\documentclass[final]{elsarticle}

\usepackage{lineno,hyperref}

\usepackage{amsthm,amssymb,epsfig,float,color}
\usepackage{siunitx}
\usepackage{upgreek}
\usepackage{geometry}
\usepackage{lipsum}
\usepackage{amsfonts}
\usepackage{booktabs} 
\usepackage{graphicx}
\usepackage{epstopdf}
\usepackage{algorithmic,multirow,mathtools}
\graphicspath{{Graphics/}}

\usepackage{float}
\usepackage{hyperref}
\usepackage{color}
\usepackage{tcolorbox}
\tcbuselibrary{xparse}
\usepackage{relsize}
\usepackage{latexsym}
\usepackage{import}
\usepackage{xifthen}
\usepackage{transparent}

\allowdisplaybreaks

\newcommand{\rmd}{\mathrm{d}}

\newcommand{\dd}[2]{\fracd{\rmd {#1}}{\rmd {#2}}}
\newcommand{\fracd}[2]{\displaystyle
{\frac{{\displaystyle{#1}}}{{\displaystyle{#2}}}}}

\newcommand{\Qe}{Q_{\mathrm{e}}}

\newcommand{\Qf}{Q_{\mathrm{f}}}

\newcommand{\Qu}{Q_{\mathrm{u}}}

\newcommand{\Xc}{X_{\mathrm{c}}}

\newcommand{\Xf}{X_{\mathrm{f}}}

\newcommand{\bC}{\boldsymbol{C}}

\newcommand{\bPhi}{\boldsymbol{\Phi}}

\newcommand{\bR}{\boldsymbol{R}}

\newcommand{\bSf}{\boldsymbol{S}_{\mathrm{f}}}

\newcommand{\bS}{\boldsymbol{S}}

\newcommand{\bzero}{\boldsymbol{0}}

\newcommand{\mumax}{\mu_{\mathrm{max}}}

\newcommand{\Qreac}{Q_{\mathrm{reac}}}

\newcommand{\sme}{\sigma_{\mathrm{e}}}

\newcommand{\vhs}{v_{\mathrm{hs}}}
\newcommand{\vrel}{v_{\mathrm{rel}}}

\newcommand{\br}{\boldsymbol{r}}
\newcommand{\bsigmaC}{\boldsymbol{\sigma}_\mathrm{\!\bC}}
\newcommand{\bsigmaS}{\boldsymbol{\sigma}_\mathrm{\!\bS}}

\definecolor{ros}{RGB}{148,35,9}   

\newcommand{\bCf}{\bC_{\rm f}}

\DeclareMathAlphabet{\mathbcal}{OMS}{cmsy}{b}{n}

\begin{document}

\begin{frontmatter}

\title{A moving-boundary model of reactive settling in wastewater treatment. Part~1: Governing equations}


\author[RBaddress]{Raimund B\"urger}
\author[JCaddress]{Julio Careaga}
\author[JCaddress]{Stefan Diehl\corref{mycorrespondingauthor}}
\cortext[mycorrespondingauthor]{Corresponding author, \texttt{stefan.diehl@math.lth.se}}
\author[RBaddress]{Romel Pineda}

\address[RBaddress]{CI${}^{\mathrm{2}}$MA and Departamento de Ingenier\'{\i}a Matem\'{a}tica, Facultad de Ciencias F\'{i}sicas y Matem\'{a}ticas, Universidad~de~Concepci\'{o}n, Casilla 160-C, Concepci\'{o}n, Chile}
\address[JCaddress]{Centre for Mathematical Sciences, Lund University, P.O.\ Box 118, S-221 00 Lund, Sweden}

\begin{abstract}   
Reactive settling is the process of sedimentation of small solid particles  in a  fluid with simultaneous reactions between the components of the solid and liquid phases. This process is  important in sequencing batch reactors (SBRs) in wastewater treatment plants. In that application the particles are biomass (bacteria; activated sludge) and the liquid contains substrates (nitrogen, phosphorus) to be removed through reactions with the biomass. The operation of an SBR in cycles of consecutive  fill, react, settle, draw, and idle stages is modelled by a system of spatially one-dimensional, nonlinear, strongly degenerate parabolic convection-diffusion-reaction  equations. This system is coupled via conditions of mass conservation to transport equations on a half line, whose origin is located at a moving boundary and that model the effluent pipe.
{An invariant-region-preserving finite difference scheme is used} to simulate operating cycles and the denitrification process within an SBR.
\end{abstract}

\begin{keyword}
convection-diffusion-reaction PDE\sep degenerate parabolic PDE\sep moving boundary\sep sedimentation\sep sequencing batch reactor
\MSC[2010] 35K65\sep 35Q35\sep 65M06\sep 76V05
\end{keyword}

\end{frontmatter}


\section{Introduction}

\subsection{Scope}

\begin{figure}[t]
\centering 
\includegraphics[width=0.99\textwidth]{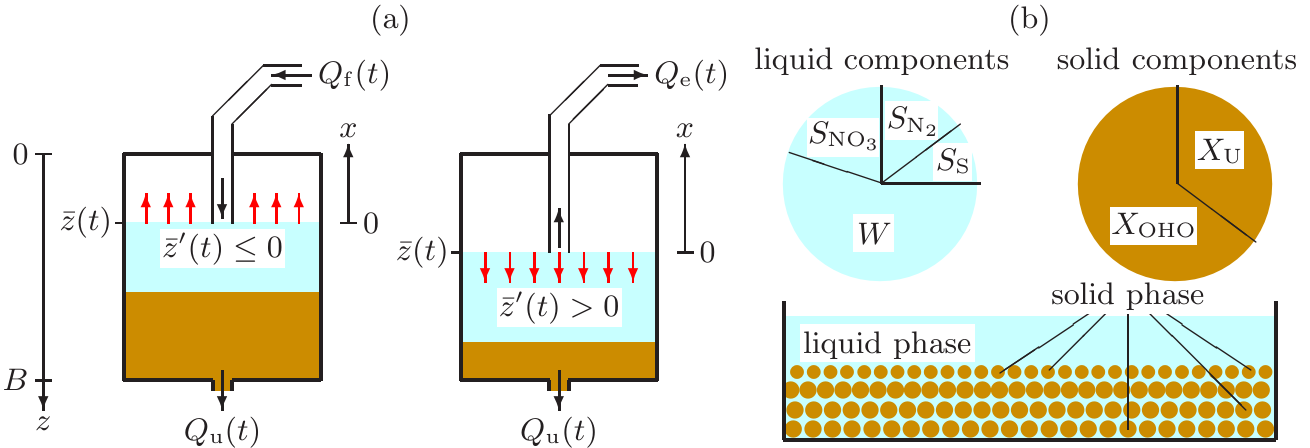}
\caption{(a) Left: Fill at the volume rate $\Qf(t)>0$ greater than the underflow rate $\Qu(t)\geq 0$ resulting in a rise of the mixture surface location~$z=\bar{z}(t)$. Right: Draw (extraction) of mixture from the surface at the rate $\Qe(t)>0$ implies a descending surface. (b) The two  phases and their  components  for the examples of denitrification in Section~\ref{sec:numex}.}\label{fig:FillDraw}
\end{figure}

\begin{figure}[t]
\centering
 \includegraphics[width=0.99\textwidth]{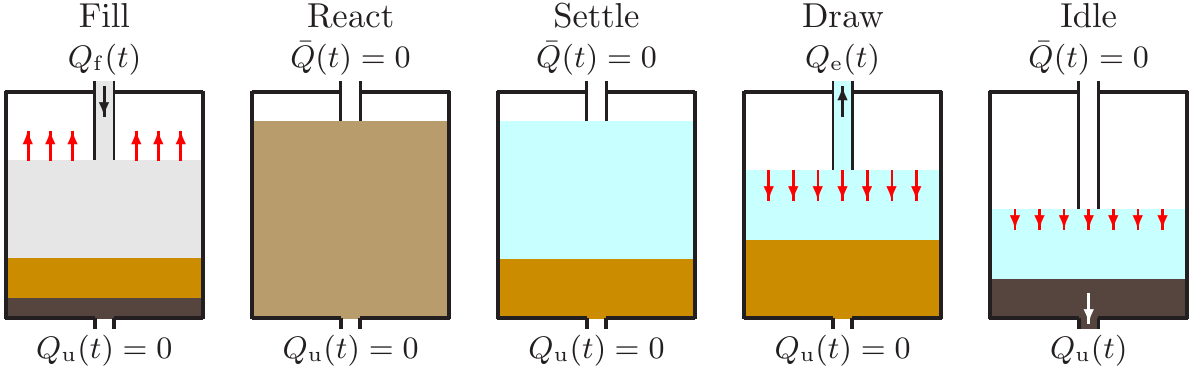}
\caption{The five stages of a cycle of an SBR.
The tank is first filled with wastewater at the volumetric flow $\Qf(t)>0$ and concentrations~$\bCf(t)$ and $\bSf(t)$.
During the react stage, biological reactions take place under complete mixing by an impeller or by aeration.
Then batch sedimentation with reactions occurs and liquid is extracted during the draw stage.
During the idle stage, some of the bottom sludge can be withdrawn and then the fill stage starts again.}
\label{fig:SBRcycle}
\end{figure}

We present a  one-dimensional model of reactive sedimentation in a tank (with a possibly varying cross-sectional area). At the bottom, the tank has a controlled outlet. 
At the surface of the mixture, a floating device allows for controlled fill or extraction of mixture; see Figure~\ref{fig:FillDraw}~(a).
The settling particles consist of several components, which react with other dissolved material components.
{The model can handle} any feed or extraction condition where the volume of mixture in the tank may vary between zero (surface at the bottom) and maximal (surface at the top).
The specific application we have in mind is a \emph{sequencing batch reactor} (SBR), which is commonly used for wastewater treatment, where batch operations of reactions and sedimentation are applied in sequence in time, with fill and draw (extraction) operations between or during these stages; see Figure~\ref{fig:SBRcycle}.
In an SBR, the particles are biomass (bacteria; activated sludge) and the dissolved materials are substrates (nitrogen, phosphorus, etc.) to be removed.
Other applications arise, for example, in mineral processing where mineral powders are flocculated by adding liquid flocculant dissolved in water.

To introduce the governing model, we let $A=A(z)$ denote the cross-sectional area of the tank that may depend on depth~$z$, where $z=0$ at the top of the tank and $z=B$ at its bottom.
The characteristic function~$\gamma$ equals one inside the mixture and zero otherwise, i.e., $\gamma(z,t)=\chi_{\{\bar{z}(t)<z<B\}}$, where $\chi_I$ is the indicator function which equals one if and only if $I$ is true, and $z=\bar{z}(t)$ is the surface location.
The unknowns are the vectors  $\smash{\bC =  \bC(z,t) =  ( C^{(1)}, \dots, C^{(k_{\boldsymbol{C}})} )^{\mathrm{T}}}$ of solid concentrations  and of $\smash{\bS = \bS(z,t)
= ( S^{(1)}, \dots, S^{(k_{\boldsymbol{S}})})^{\mathrm{T}} }$ concentrations of soluble components. These vectors make up the components of the 
solid and liquid phase, respectively. 
With $X=X(z,t) = \smash{C^{(1)} (z,t) + \dots + C^{(k_{\boldsymbol{C}})} }(z,t)$, the tank can be 
modelled as the following system of convection-diffusion-reaction equations, 
where $t>0$ is time: 
\begin{align} \label{eq:model} 
	\begin{split} 
&A(z)\partial_t{\bC} + \partial_z\big(A(z) \mathcal{F}_{\boldsymbol{C}} (X,z,t)  \bC\big)- \partial_{z}\big(A(z)\gamma(z,t)(\partial_z{D(X)}) \boldsymbol{C} \big)   \\
&
\qquad\qquad   =  \delta \bigl(z-\bar{z} (t) \bigr) \Qf(t)\bCf(t) + \gamma(z,t)A(z)\bR_{\bC}(\bC,\bS), 
\\
& A(z)\partial_t{\bS} + \partial_{z}\big( A(z) \mathcal{F}_{\boldsymbol{S}} (X,z,t)  \bS\big)   \\
& \qquad\qquad  =  \delta \bigl(z-\bar{z} (t) \bigr) \Qf(t)\bCf(t) + \gamma(z,t)A(z)\bR_{\bS}(\bC,\bS).
  \end{split}
\end{align}
This system is coupled to a model of the effluent pipe consisting of convective transport equations on a half line~$x\geq0$, where~$x=0$ is attached to the moving boundary~$z=\bar{z}(t)$ (cf.\ Figure~\ref{fig:FillDraw}~(a)). 
The coupling conditions between the systems are mass-preserving algebraic equations with fluxes on the~$z$- and $x$-axes.
The scalar functions $\mathcal{F}_{\boldsymbol{C}}$ and $\mathcal{F}_{\boldsymbol{S}}$ in~\eqref{eq:model} depend 
   nonlinearly on~$X$ and represent portions of the solid and
     liquid phase velocity, respectively.  The scalar 
    function $D$ models sediment compressibility. 
The terms with the delta function $\delta(z-\bar{z}(t))$ model the operation of the  device  floating on the surface during feed of mixture. 
The last term of each {partial differential equation}  (PDE) contains the reaction rates (local increase of mass per unit time and volume) $\bR_{\bC}$ and $\bR_{\bS}$.
The full PDE model is specified in Section~\ref{sec:model}.
During the react stage of an SBR (see Figure~\ref{fig:SBRcycle}), full mixing occurs and the system of PDEs~\eqref{eq:model} reduces to a system of ordinary differential equations (ODEs).

{The purpose of this work is to derive the model and present simulations when the reaction terms model a simplified denitrification process in wastewater treatment, which occurs when no oxygen is present.}
The main difficulties for the analysis of the {entire PDE} model arise partly from the presence of a moving boundary where both a source is located and a half-line model attached, and partly from strong type degeneracy; the function $D$ is zero for $X$-values on an interval of positive length.
{The simulations are made with a new positivity-preserving numerical scheme that handles these difficulties and is presented in~\cite{bcdp_part2}.}

\subsection{Related work}

The SBR technology has been used for hundred years and been a topic {of intense} research~\cite{Singh2010}.
Its usage for wastewater treatment {is outlined}  in many handbooks (e.g., \cite{chen2020book,droste,metcalf}). 
Furthermore, it is also employed for recovery of selenium~\cite{Song2021}, radioactively labelled pharmaceuticals~\cite{Popple2016}, nitrogen removal processes~\cite{Ni2014}, pharmaceutically active compounds~\cite{Wang2011}, synthetic chemical components~\cite{Hu2005}, swine manure slurry~\cite{Mass2000}, applications in the petrochemical industry~\cite{Caluw2016}, and saline wastewater treatment~\cite{Amin2014} among others.

Most treatments in the literature on mathematical models related to SBRs focus  on ODEs modelling the reactions   by established activated sludge models~\cite{Freytez2019,Henze2000ASMbook,Kauder2007,Meadows2019},  
optimization and control problems~\cite{Gajardo2009,Kim2009,Piotrowski2019,Souza2008,Pambrun2008},  and statistical methods~\cite{Kocijan2013,Li2012}.
Less consideration has been laid on the sedimentation in an SBR during which reactions occur.
Models of  reactive settling in continuously operated secondary settling tanks (SSTs) based on PDEs are presented in~\cite{bcdimamat21,SDcace_reactive,SDm2an_reactive}
 (see also references cited in these works). 
It is worth pointing out that the SBR model differs from an SST model. 
In an SBR, $\Qu(t)$, $\Qf(t)$ and $\Qe(t)$ are given independent control functions giving rise to a moving surface, whereas in an SST, only two of these are known and the third, often~$\Qe(t)$, is defined by the other two and possibly by volume-changing reactions in the tank~\cite{SDwatres3}.

\subsection{Outline of the paper}

In Section~\ref{sec:model}, we derive the model. 
Section~\ref{sec:applSBR} contains a derivation of the ODE model for full mixing.
In Section~\ref{sec:numex}, we show two numerical examples of SBR operation with a constant cross-sectional area (cylindrical vessel)
and a variable.
Some conclusions are collected in Section~\ref{sec:concl}.

\section{Derivation of the model of reactive settling}\label{sec:model}

\subsection{Preliminaries} \label{subsec:notations}

The solid phase consists of flocculated particles of $k_{\boldsymbol{C}}$ types   with concentrations~$\smash{C^{(1)} ,\dots, C^{(k_{\boldsymbol{C}})}}$.
 The components of the liquid phase are water of concentration~$W$ and $k_{\bS}$ dissolved substrates of concentrations $\smash{S^{(1)} ,\dots, S^{(k_{\bS})}}$ (cf.\  Figure~\ref{fig:FillDraw}~(b)). 
The total concentrations of solids~$X$ and liquid~$L$ are
\begin{align} \label{eq:XLdef} 
 X :=  C^{(1)} + \dots  + C^{(k_{\boldsymbol{C}})} ,\quad L := 
  W + S^{(1)} + \dots + S^{(k_{\boldsymbol{S}})}.
\end{align} 
All these concentrations depend on~$z$ and~$t$, and our notation is the same as in~\cite{bcdimamat21}.

For computational purposes, we define a maximal concentration~$\hat{X}$ of solids and  assume that the density of all solids is the same, namely~$\rho_X>\hat{X}$.
Similarly, we assume that the liquid phase has density~$\rho_L<\rho_X$, typically the density of water.
{The reaction terms for all components are collected in the vectors 
\begin{equation*}
\bR_{\bC}(\bC,\bS)=\bsigmaC\br(\bC,\bS),\qquad
\bR_{\bS}(\bC,\bS)=\bsigmaS\br(\bC,\bS),
\end{equation*}
which model the increase of solid and soluble components, respectively, where $\bsigmaC$ and $\bsigmaS$ are constant stoichiometric matrices and $\br(\bC,\bS)\geq\bzero$ is a vector of non-negative reaction rates, which are assumed to be bounded and Lipschitz continuous functions.}
We set 
\begin{align*} 
\tilde{R}_{\bC}(\bC,\bS) :=  R_{\bC}^{(1)}(\bC,\bS)
+\dots + R_{\bC}^{(k_{\bC})}(\bC,\bS)
\end{align*} 
 (analogously for~$\smash{\tilde{R}_{\bS}(\bC,\bS)}$).
The water concentration $W$ is not active in any reaction.

We let $v_X$ and $v_L$ denote the velocities of the  solid and liquid phases, respectively. It is assumed that the relative velocity
$v_X-v_L=:\vrel=\vrel(X,\partial_z X,z)$ 
is given by a constitutive function of~$X$ and its spatial derivative~$\partial_z X$, modelling hindered and compressive settling; see Section~\ref{subsec:constitutive}.
{The reason for the dependence on the total concentration $X$ (and its spatial derivative) is that the particles are flocculated consisting of all solid components.
All components within a particle settle with the same velocity.}
{To obtain non-negative concentrations of the vector~$\bC$ (analogously for~$\bS$), we let
\begin{equation*}
I_{\bC,k}^-:=\big\{l\in\mathbb{N}:\sigma_{\bC}^{(k,l)}<0\big\}
\end{equation*}
denote the set of indices~$l$ that have negative stoichiometric coefficients  and assume the following:
\begin{equation}\label{eq:positivity}
\text{if $l\in I_{\bC,k}^-$, then $r^{(l)}(\bC,\bS)=\bar{r}^{(l)}(\bC,\bS)C^{(k)}$ with $\bar{r}^{(l)}$ bounded.}
\end{equation}
The assumption~\eqref{eq:positivity} is natural, since it implies that (analogously for~$\bS$)
\begin{equation*} 
R^{(k)}_{\bC}(\bC,\bS)\bigr|_{C^{(k)}=0} \geq 0 \quad \text{for $k=1, \dots, k_{\boldsymbol{C}}$,}\label{eq:assumptionRC}
\end{equation*}
which means that the system of ODEs 
\begin{equation*}
\frac{\rmd}{\rmd t} \begin{pmatrix}
\bC\\ \bS
\end{pmatrix} = \begin{pmatrix}
\bR_{\bC}(\bC,\bS)\\ \bR_{\bS}(\bC,\bS)
\end{pmatrix}
\end{equation*}
has non-negative solutions if the initial data are non-negative~\cite{Formaggia2011}.}

It is assumed that in the inlet and outlet pipes no reactions take place and all the components have the same velocity.
At the bottom, $z=B$, one can withdraw mixture at a given volume rate  $\Qu(t)\geq 0$.
The underflow region $z>B$ is for simplicity modelled by setting $A(z):= A(B)$, since we are only interested in the underflow concentration $\bC_\mathrm{u}(t)$, which is an outcome of the model (analogously for~$\bS_\mathrm{u}(t)$).

At the surface of the mixture, $z=\bar{z}(t)$, we model a floating device connected to a pipe through which one can feed the tank with given volume rate~$\Qf(t)$ and  feed concentrations~$\bCf(t)$ and~$\bSf(t)$; see Figure~\ref{fig:FillDraw}.
This gives rise to a source term in the model equation with the fluxes $\Qf(t)\bCf(t)$ and $\Qf(t)\bSf(t)$.
Alternatively, this floating device allows one to extract mixture at a given volume rate $\Qe(t)>0$ through the same pipe; hence, one cannot fill and extract simultaneously.
If~$[0,T]$ denotes the total time interval of modelling (and simulation in Section~\ref{sec:numex}), we assume that $T:={T}_\mathrm{e}\cup{T}_\mathrm{f}$, where
\begin{align*}
{T}_\mathrm{e}&:= \bigl\{t\in\mathbb{R}_+:\Qe(t)>0,\Qf(t)=0 \bigr\},
\qquad
{T}_\mathrm{f}:= \bigl\{t\in\mathbb{R}_+:\Qe(t)=0,\Qf(t)\geq0 \bigr\}.
\end{align*}
{Periods when there is neither extraction nor filling are for convenience included in ${T}_\mathrm{f}$.}
When~$t\in{T}_\mathrm{e}$, we model the extraction flow in the effluent pipe by a moving coordinate system; a half line~$x\geq 0$, where $x=0$ is attached to $z=\bar{z}(t)$.
Along this half line, we denote the solids concentration by~$\boldsymbol{\tilde{C}}=\boldsymbol{\tilde{C}}(x,t)$.
The effluent concentration~$\bC_\mathrm{e}(t):=\boldsymbol{\tilde{C}}(0^+,t)\chi_{\{t\in{T}_\mathrm{e}\}}$ is also a model outcome (analogously for~$\bS_\mathrm{e}(t)$).

It is convenient to define  the volume fractions 
\begin{equation} \label{eq:phidef}
\phi:= X/\rho_X ,\qquad \phi_L:= L/ \rho_L,\qquad \phi_\mathrm{M}:=\phi+\phi_L,
\end{equation} 
where the volume fraction of the mixture
satisfies $\phi_\mathrm{M}=\chi_{\{z>\bar{z}(t)\}}$.
Below the surface, $\phi+\phi_L=1$, or  equivalently, $L=\rho_L(1-X/\rho_X)$.
The same holds for the feed concentrations.
For known~$\bC$ and~$\bS$, \eqref{eq:XLdef} implies   the water concentration  
\begin{equation}  \label{eq:W}
W = \rho_L(1-X/\rho_X)-\big(S^{(1)}+\cdots+S^{(k_{\bS})}\big).
\end{equation} 
This concentration is not part of any reaction and can be computed afterwards.

The volume of the mixture is defined by 
\begin{align} \label{eq:V}
\bar{V}(t):= V \bigl(\bar{z}(t)\bigr), \quad \text{where} \quad 
V(z):= \int_{z}^{B}A(\xi)\,\mathrm{d}\xi \quad \text{for $ 0\leq z\leq B$.} 
\end{align} 
The function~$V$  is invertible since $V'(z)=-A(z)<0$; in particular,
\begin{align} \label{eq:VA}
\bar{V}'(t)=V' \bigl(\bar{z}(t) \bigr)\bar{z}'(t)=-A \bigl(\bar{z}(t) \bigr)\bar{z}'(t).
\end{align}

\subsection{Balance laws} \label{subsec:pdemodel}

The balance laws for all components in local form imply
\begin{subequations}\label{eq:mod}
\begin{align}
\partial_{t} \bigl(A(z)\bC \bigr)+\partial_{z}\big(A(z)v_X\bC\big) & =
 \delta \bigl(z-\bar{z}(t) \bigr)\Qf\bCf + \gamma(z,t)A(z)\bR_{\bC},  \quad z\in\mathbb{R},\label{eq:govC} \\
\partial_{t} \bigl(A(z)\bS \bigr)+\partial_{z}\big(A(z)v_L\bS\big) & = \delta \bigl(z-\bar{z}(t) \bigr)\Qf\bSf + \gamma(z,t)A(z)\bR_{\bS},  \quad z\in\mathbb{R}.\label{eq:govS}
\end{align}
\end{subequations}
This system along with~$v_{\mathrm{rel}} = v_X - v_L$ and \eqref{eq:phidef} are $k_{\bC}+k_{\bS}+2$ equations for the same number of scalar unknowns,  i.e.,  the components of $\bC$ and $\bS$, plus $v_X$ and $v_L$.
It is coupled to the following model of the effluent pipe with cross - area~$A_\mathrm{e}$:
\begin{subequations}\label{eq:mod_e}
\begin{alignat}2
A_\mathrm{e}\partial_{t} \boldsymbol{\tilde{C}} + \Qe\partial_{x}\boldsymbol{\tilde{C}} &= \bzero, &\quad& x> 0, \label{eq:govCpipe}\\
A_\mathrm{e}\partial_{t} \boldsymbol{\tilde{S}} + \Qe\partial_{x}\boldsymbol{\tilde{S}} &= \bzero, &\quad& x> 0, \label{eq:govSpipe}\\
-\Qe(t)\boldsymbol{\tilde{C}}(0^+,t) &= A\bigl(\bar{z}(t)^+ \bigr)\big(
v_X |_{z=\bar{z}(t)^+} -\bar{z}'(t)\big)\bC \bigl(\bar{z}(t)^+,t \bigr),&\quad&\label{eq:connC}\\
-\Qe(t)\boldsymbol{\tilde{S}}(0^+,t) &= A \bigl(\bar{z}(t)^+ \bigr)\big(
 v_L |_{z=\bar{z}(t)^+} -\bar{z}'(t)\big)\bS\bigl(\bar{z}(t)^+,t\bigr).&\quad&\label{eq:connS}
\end{alignat}
\end{subequations}
The coupling equations~\eqref{eq:connC} and \eqref{eq:connS} preserve  mass at the surface during extraction periods.
The purpose of~\eqref{eq:mod_e} is to define the concentrations during periods of extraction when $\Qe(t)>0$.
For~$t>0$ the outlet concentrations are given by
\begin{alignat}2
\bC_{\rm u}(t)&:=\bC(B^+,t),&\qquad
\bS_{\rm u}(t)&:=\bS(B^+,t), \label{eq:u_conc_def}\\
\bC_\mathrm{e}(t)&:=\boldsymbol{\tilde{C}}(0^+,t) \chi_{\{t\in{T}_\mathrm{e}\}},&\qquad
\bS_\mathrm{e}(t)&:=\boldsymbol{\tilde{S}}(0^+,t) \chi_{\{t\in{T}_\mathrm{e}\}}.\label{eq:e_conc_def}
\end{alignat}
The transport PDEs~\eqref{eq:govCpipe} and \eqref{eq:govSpipe} are easily solved once the boundary data~\eqref{eq:e_conc_def} are known, which in turn have to satisfy~\eqref{eq:connC} and \eqref{eq:connS}.
The right-hand sides of the latter equations are however nonlinear functions of~$\bC(\bar{z}(t)^+,t)$ (via $v_X$ and $v_L$) and to obtain unique boundary concentrations on either side of a spatial discontinuity, an additional entropy condition is needed.
Our experience is, however, that correct concentrations can be obtained by a conservative and monotone numerical method~\cite{Burger&K&T2005a}.

The volume-average bulk velocity is defined by $q(z,t):= (\phi v_X+\phi_Lv_L)\chi_{\{z>\bar{z}(t)\}}$.
Since $\phi_\mathrm{M}=0$ for $z<\bar{z}(t)$, we have $\phi=\phi_L=0$ there.
Summing all the equations of~\eqref{eq:govC} and \eqref{eq:govS}, respectively, and using~\eqref{eq:XLdef} and \eqref{eq:phidef}, we get the scalar PDEs
\begin{align*} 
\partial_{t} \bigl(A(z)\rho_X\phi \bigr)+\partial_{z} \bigl(A(z)\rho_X\phi v_X \bigr)     &= \delta 
 \bigl(z-\bar{z}(t) \bigr)\rho_X{\phi_\mathrm{f}}\Qf+\gamma(z,t)A(z)\tilde{R}_{\bC},\\
\partial_{t} \bigl(A(z)\rho_L\phi_L \bigr)+\partial_{z} \bigl(A(z)\rho_L\phi_L v_L \bigr) &= \delta
 \bigl(z-\bar{z}(t) \bigr)\rho_L\phi_{L,\mathrm{f}}\Qf+\gamma(z,t)A(z)\tilde{R}_{\bS}.
\end{align*}
Dividing away the constant densities and adding 
the results, we get  
\begin{equation} \label{eq:temp1}
\partial_{t} \bigl(A(z)\phi_\mathrm{M} \bigr) + \partial_{z} \bigl(A(z)q \bigr)= \delta
 \bigl(z-\bar{z}(t) \bigr)\phi_{\mathrm{M,f}}\Qf+\gamma(z,t)A(z)\mathcal{R},
\end{equation}
where $\mathcal{R} :=  \tilde{R}_{\bC} / \rho_X +  \tilde{R}_{\bS} / \rho_L$, and $\phi_{\mathrm{M,f}}=1$ by definition.
The first term of~\eqref{eq:temp1} is
\begin{equation*}
\partial_{t} \bigl(A(z)\phi_\mathrm{M} \bigr) = A(z)\partial_{t}\chi_{\{z>\bar{z}(t)\}} 
= -A(z) \delta \bigl(z-\bar{z}(t) \bigr)\bar{z}'(t).
\end{equation*}
The same procedure for the algebraic equations~\eqref{eq:connC} and~\eqref{eq:connS} ($t\in I_\mathrm{e}$)
 yields  
\begin{align} \label{eq:tempalg}
-\Qe(t)=A\big(\bar{z}(t)^+\big)\big(q(\bar{z}(t)^+,t) - \bar{z}'(t)\big).
\end{align} 
Integrating~\eqref{eq:temp1} (with or without the source term) from $z\in(\bar{z}(t)^+,B)$ to $B$, we get
\begin{align*}
A(B)q(B,t) - A(z)q(z,t) = \int_{z}^{B}A(\xi) \mathcal{R}\big(\bC(\xi,t),\bS(\xi,t)\big)\,\rmd\xi
=: \Qreac(z,t;\bC,\bS),
\end{align*}
where $A(B)q(B,t) = \Qu(t)$. Hence, inside the mixture, i.e., in the interval $(\bar{z}(t),B)$, the volume-average velocity~$q$ is given by
$A(z)q(z,t)  = \Qu(t) - \Qreac(z,t;\bC,\bS)$. 
In view of this equation and $q= (\phi v_X+\phi_Lv_L)\chi_{\{z>\bar{z}(t)\}}$, we integrate~\eqref{eq:temp1} from $z=\bar{z}(t)-h$ to $\bar{z}(t)+h$ and let $0<h\rightarrow 0$ to get
\begin{align} \label{eq:temp2}
-A \bigl(\bar{z}(t) \bigr)\bar{z}'(t) + \Qu(t) - \left.\Qreac(z,t;\bC,\bS)\right|_{z=\bar{z}(t)} = \Qf(t),
\end{align} 
where the first term can be written $\bar{V}'(t)$; see \eqref{eq:VA}.
For $t\in {T}_\mathrm{e}$, \eqref{eq:tempalg} implies
\begin{align} \label{eq:temp3}
-\Qe(t)=\Qu(t)-\left.\Qreac(z,t;\bC,\bS)\right|_{z=\bar{z}(t)} +V'(t).
\end{align} 

The term $\Qreac$ seems to be negligible in the application to wastewater treatment~\cite{SDm2an_reactive} and we set $\Qreac:=0$ from now on.
Consequently, we define
\begin{equation} \label{eq:qdef2}
q(z,t) :=\bigl( \Qu(t) / A(z) \bigr) \chi_{\{z>\bar{z}(t)\}}.
\end{equation} 
Then \eqref{eq:temp2} and \eqref{eq:temp3} can be written (with $\Qreac=0$) as 
\begin{align} \label{eq:Vprime}
\bar{V}'(t) = \bar{Q}(t) - \Qu(t),\quad\text{where}\quad
\bar{Q}(t)=\begin{cases}
-\Qe(t)<0 & \text{if $t\in {T}_\mathrm{e}$,} \\
\Qf(t)\geq 0 & \text{if $t\in {T}_\mathrm{f}$.} 
\end{cases}
\end{align} 
Solving this ODE and utilizing~\eqref{eq:V}, we obtain
\begin{align} \label{eq:zbar}
\bar{z}(t)=V^{-1}\left(
\bar{V}(0)+ \int_{0}^{t}\big(\bar{Q}(s) -\Qu(s)\big)\,\mathrm{d}s\right).
\end{align} 
Alternatively, $\bar{z}(t)$ can be obtained from (see~\eqref{eq:VA})
\begin{align} \label{eq:zbarprime}
\bar{z}'(t)=\frac{\Qu(t)-\bar{Q}(t)}{A(\bar{z}(t))}.
\end{align} 

\subsection{Constitutive functions for hindered and compressive settling}\label{subsec:constitutive}

The  surface location $z=\bar{z}(t)$ is now  specified,
so we may focus on the mixture in $z>\bar{z}(t)$.
For the given functions~$q$ and the relative velocity~$\vrel$, we set $v:=(1-\phi)\vrel$, where $\phi=X/\rho_X$, and obtain from~$v_{\mathrm{rel}} = v_X - v_L$ and $q= (\phi v_X+\phi_L v_L)\chi_{\{z>\bar{z}(t)\}}$ the  phase velocities 
\begin{equation} 
v_X = q+(1-\phi)\vrel = q +v \quad \text{and}
\quad   v_L = q-\phi\vrel = q - \dfrac{\phi}{1-\phi}v   \label{eq:vX}
\end{equation} 
of the solid and fluid, respectively. 
We assume that the relative velocity $\vrel=v/(1-\phi)$ {is given through the following commonly used expression~\cite{SDwatres3,Burger&K&T2005a} for~$v$}: 
\begin{equation}\label{eq:v}
v(X,\partial_z X,z,t) 
:= \gamma(z,t)\vhs(X)\left(
1-\dfrac{\rho_X\sme'(X)}{Xg\Delta\rho}\partial_z{X}
\right)
= \gamma(z,t)\big(\vhs(X) - \partial_z{D(X)}\big),
\end{equation}
where
\begin{equation*}
D(X) := \int_{X_c}^{X}d(s)\,{\rm d}s,
\qquad
d(X):=\vhs(X)\dfrac{\rho_X\sme'(X)}{gX\Delta\rho}. 
\end{equation*}
Here $\Delta\rho:=\rho_X-\rho_L$, $g$ is the acceleration of gravity, $\vhs=\vhs(X)$ is the hindered-settling velocity, which is assumed to be decreasing and satisfy $\vhs(\hat{X})=0$. {Moreover,} 
$\sme=\sme(X)$ is the effective solids stress, which satisfies 
\begin{align*} 
\sme'(X) \begin{cases} = 0 & \text{for $X \leq \Xc$,} \\
    > 0 & \text{for $X>\Xc$,} \end{cases} \end{align*}   
where $\Xc$ is a critical concentration above which the particles touch each other and form a network that can bear a certain stress.
Note that $d(X) =0$ for $X \leq X_{\mathrm{c}}$, which causes the strongly degenerate type behaviour. 

{The determination of these constitutive assumptions including the critical concentration~$X_\mathrm{c}$ is a topic in itself; see e.g.\ \cite{SDAPNUM1,SDwatres_Torfs2017} and references therein.}

\subsection{Model equations in final form} \label{subsec:gov}

With~$q$ defined by~\eqref{eq:qdef2}, we define and use~\eqref{eq:vX} and~\eqref{eq:v} to write the velocities
\begin{align} \label{eq:modelterms} 
\begin{split} 
 \mathcal{F}_{\boldsymbol{C}} (X,z,t)  & :=  q(z,t) + \gamma(z,t) \vhs(X),  \\
 \mathcal{F}_{\boldsymbol{S}} (X,z,t) & := \frac{ \rho_X q(z,t) - ( q(z,t) + \gamma(z,t) \vhs(X) )  X }{\rho_X - X},
\end{split}
\end{align}
and then express the total mass fluxes of the balance laws~\eqref{eq:mod} in light of~\eqref{eq:vX}:
\begin{equation*} 
\begin{split} 
\bPhi_{\bC}  := \bPhi_{\bC}(\bC,X,\partial_z X,z,t) &:= A(z)v_{X}(X,\partial_z X,z,t)\bC\\
& = A(z)\big(\mathcal{F}_{\boldsymbol{C}}(X,z,t) - \gamma(z,t)\partial_z{D(X)}\big), \\
\bPhi_{\bS}  := \bPhi_{\bS}(\bS,X,\partial_z X,z,t) &:= A(z)\frac{\rho_Xq-v_X X}{\rho_X-X}\bS = A(z)\mathcal{F}_{\boldsymbol{S}}(X,z,t). 
\end{split} 
\end{equation*} 
Then we define and rewrite the right-hand side of~\eqref{eq:connC} with~\eqref{eq:zbarprime}, \eqref{eq:vX} and \eqref{eq:v}:
\begin{equation}\label{eq:PhiC_e}
\begin{aligned} 
\bPhi_{\bC,\mathrm{e}}(z,t)
&:= 
A(\bar{z}(t))\big(v_{X}(\bar{z}(t)^+,t) -\bar{z}'(t)\big)\bC(\bar{z}(t)^+,t)\\
&= \left.
\Big(A(z)\big(q+\vhs(X) - \partial_z{D(X)}\big) 
- \Qu-\Qe\Big)\bC
\right|_{z=\bar{z}(t)^+}\\
&= \left.\Big(
A(z)\big(
\vhs(X) - \partial_z{D(X)}\big)
- \Qe
\Big)\bC
\right|_{z=\bar{z}(t)^+}.
\end{aligned}
\end{equation} 
Analogously, we define $\bPhi_{\bS,\mathrm{e}}$ corresponding to~\eqref{eq:connS}:
\begin{equation}\label{eq:PhiS_e}
\begin{aligned} 
\bPhi_{\bS,\mathrm{e}}(z,t)
&:= 
A(\bar{z}(t))\big(v_{L}(\bar{z}(t)^+,t) -\bar{z}'(t)\big)\bS(\bar{z}(t)^+,t)\\
&= \left.
\left(A(z)\left(
q - \frac{X(\vhs(X) - \partial_z{D(X)})}{\rho_X-X}
\right) - (\Qu+\Qe)\right)\bS
\right|_{z=\bar{z}(t)^+}\\
&= \left.
- \left(A(z)
 \frac{X(\vhs(X) - \partial_z{D(X)})}{\rho_X-X}
+ \Qe\right)\bS
\right|_{z=\bar{z}(t)^+}.
\end{aligned}
\end{equation} 

The final model can now be described as follows.
Given the in- and outgoing volumetric flows, one computes the surface level by~\eqref{eq:zbar} or~\eqref{eq:zbarprime}.
The concentrations~$\bC$ and $\bS$ are given by the system~\eqref{eq:mod}, which can be written as
\begin{subequations}\label{finalmod}
\begin{align}
A(z)\partial_t{\bC}+\partial_z{\bPhi_{\bC}} & = \delta \bigl(z-\bar{z}(t)\bigr)\Qf\bCf + \gamma(z,t)A(z)\bR_{\bC},  \label{finalmod_a} \\
A(z)\partial_t{\bS} +\partial_z{\bPhi_{\bS}} & = \delta \bigl(z-\bar{z}(t)\bigr)\Qf\bSf
+ \gamma(z,t)A(z)\bR_{\bS}, \label{finalmod_b}
\end{align} 
\end{subequations}
or as \eqref{eq:model}.
The water concentration $W$ can always be calculated from~\eqref{eq:W}.
Note that~$W$ is not present in \eqref{eq:model} or  \eqref{finalmod}.
The effluent and underflow concentrations are given by~\eqref{eq:u_conc_def} and \eqref{eq:e_conc_def}, respectively.
No initial data are needed for the outlet concentrations, but for the following: 
\begin{equation*}
\bC^0= \big(C^{(1),0},C^{(2),0},\dots, C^{(k_{\bC}),0}\big)^{\mathrm{T}},\qquad
\bS^0= \big(S^{(1),0},S^{(2),0},\dots, S^{(k_{\bS}),0}\big)^{\mathrm{T}}.
\end{equation*}

\section{Application to sequencing batch reactors}\label{sec:applSBR}

An SBR cycle consists of five stages; see Figure~\ref{fig:SBRcycle}.
During some of these periods, mixing may occur due to aeration or the movement of an impeller.
For the sake of simplicity, we ignore partial mixing and exemplify the cases of either no mixing or full mixing.
The PDE model~\eqref{finalmod} {excludes} mixing and we next derive the special case of full mixing.

\subsection{Model during a full mixing react stage}

Full mixing means that the relative velocity~$v_\mathrm{rel}$ is negligible.
We set $v_\mathrm{rel}\equiv 0$ and assume that concentrations only depend on~$t$ (below the surface).
Then $v_X=v_L=q$, hence,
\begin{equation} \label{eq:react_fluxes}
\bPhi_{\bC}=A(z)q(z,t)\bC=\Qu(t)\chi_{\{z>\bar{z}(t)\}}\bC,\qquad
\bPhi_{\bS}=\Qu(t)\chi_{\{z>\bar{z}(t)\}}\bS.
\end{equation} 
Integrating the PDEs~\eqref{finalmod} from $\bar{z}(t)^-$ to $B$, one gets the governing ODEs.
The integral of the time-derivative term of~\eqref{finalmod_a} can, by means of~\eqref{eq:zbarprime}, be written as 
\begin{align*}
\int_{\bar{z}(t)}^{B}\frac{\rmd\bC}{\rmd t}A(\xi)\,\rmd\xi
& =
\frac{\rmd}{\rmd t}\int_{\bar{z}(t)}^{B}\bC(t)A(\xi)\,\rmd\xi 
- \big(-\bC(t)A(\bar{z}(t))\bar{z}'(t)\big)\\
& =
\frac{\rmd\bC(t)}{\rmd t}\bar{V}(t) 
- \bC(t)\big(\Qu(t)-\bar{Q}(t)\big).
\end{align*}
The spatial-derivative term of~\eqref{finalmod_a} becomes, with \eqref{eq:react_fluxes},
\begin{equation*}
\int_{\bar{z}(t)}^{B}\frac{\mathrm{d} \bPhi_{\bC}}{\rmd z}\,\rmd\xi
=\Qu(t)\bC(t)-\bzero.
\end{equation*}
The same can be done for the substrate equations and we obtain the following system of ODEs for the homogeneous concentrations in $\bar{z}(t)<z<B$:
\begin{subequations}\label{eq:ODEmodel}
\begin{align}
\bar{V}(t)\dd{\bC}{t} & = \big(\Qu(t)-\bar{Q}(t)\big)\bC + 
\Qf(t)\bCf(t) + \bar{V}(t)\bR_{\bC},  \label{eq:ODEmodelC} \\
\bar{V}(t)\dd{\bS}{t} & = \big(\Qu(t)-\bar{Q}(t)\big)\bS + 
\Qf(t)\bSf(t) + \bar{V}(t)\bR_{\bS},  \label{eq:ODEmodelS}
\end{align}
\end{subequations}
where all the concentrations depend only on time, since they are averages (below the surface).
As before, $W$ can be obtained afterwards from~\eqref{eq:W}.
In the region $0<z<\bar{z}(t)$ all concentrations are zero.
Because of~\eqref{eq:connC} and \eqref{eq:PhiC_e}, we have  $\bC_{\rm u}(t)=\bC(t)$ and $\bC_{\rm e}(t)=\bC(t)\chi_{\{t\in{T}_\mathrm{e}\}}$ (analogously for $\bS$).
The system~\eqref{eq:ODEmodel} thus models a completely stirred tank with reactions, possibly a moving upper boundary because of in- and outflow streams.

\section{Simulations} \label{sec:numex}

{We use the novel numerical method in~\cite{bcdp_part2}.}
To exemplify the entire SBR process, we use the same model for denitrification as in~\cite{SDcace_reactive} with two solid components: ordinary heterotrophic organisms $X_{\rm OHO}$ and undegradable organics $X_{\rm U}$; and three soluble components: nitrate $S_{\rm NO_3}$, readily biodegradable substrate $S_{\rm S}$ and nitrogen $S_{\rm N_2}$ {(where we identify the denomination of a component and the corresponding concentration variable)}. Thus, we utilize 
$\bC = (
             X_{\rm OHO}, 
             X_{\rm U})^{\mathrm{T}}$ and 
             $\bS  =  (
             S_{\rm NO_3}, 
             S_{\rm S} , 
             S_{\rm N_2} 
              )^{\mathrm{T}} $, corresponding to 
              $k_{\boldsymbol{C}} =2$ and
  $k_{\boldsymbol{S}} =3$, respectively.  
Shortly described, the denitrification process converts nitrate ($\rm NO_3$) to nitrogen gas ($\rm N_2$) by   a series of reactions involving the particulate biomass.
{Since denitrification occurs without the presence of oxygen, the mixing during the react stage of the SBR process is achieved by an impeller.
The reaction vectors $\bR_{\bC}=\bsigmaC\br$ and $\bR_{\bS}=\bsigmaS\br$ have the stoichiometric matrices and reaction-rate vector
\begin{equation*}
\bsigmaC  = 
\begin{bmatrix}
1&-1\\0&f_{\rm P}
\end{bmatrix},\qquad
\bsigmaS  =
\begin{bmatrix}
-\bar{Y}&0\\
-1/Y&1-f_{\rm P}\\
\bar{Y}&0
\end{bmatrix},\qquad
\br=X_{\rm OHO}
\begin{pmatrix}
\mu(\bS)\\
b
\end{pmatrix}.
\end{equation*}
}Here $\bar{Y} = {(1-Y)}/{(2.86Y)}$ where $Y = 0.67$ is a yield factor, $b = 6.94\times 10^{-6}\,\mathrm{s}^{-1}$ is the decay rate of heterotrophic organisms, and $f_\mathrm{P} = 0.2$ is the portion of these that decays to undegradable organics.
The growth rate function
\begin{equation*}
 \mu(\bS) = \mumax \dfrac{S_{\rm NO_3}}{K_{\rm NO_3}+S_{\rm NO_3}} \dfrac{S_{\rm S}}{K_{\rm S}+S_{\rm S}}
\end{equation*}
has the parameters $\mumax = 5.56\times 10^{-5}\,\mathrm{s}^{-1}$, $K_{\rm NO_3} = 5\times 10^{-4}\, \rm kg/m^3$ and $K_{\rm S} = 0.02\, \rm kg/m^3$. 
{The first component of $\bR_{\bC}$ models partly the growth of heterotrophic organisms ($X_\mathrm{OHO}$) due to consumption of substrates with a rate proportional to its concentration with proportionality coefficient $\mu(\bS)$ and partly the decay the with coefficient~$b$.
The consumption of substrates ($S_{\rm NO_2}, S_\mathrm{S}$) is modelled by negative terms proportional to $\mu(\bS)$.}

The maximal total solids concentration is set to $\hat{X}=30\,\mathrm{kg}/\mathrm{m}^3$, a value our simulated solutions never reach.
The constitutive functions used in all simulations are 
\begin{equation*}
\vhs(X) := \frac{v_0}{ 1 + (X/ \breve{X})^\eta}, \qquad
\sigma_{\mathrm{e}}(X):=\alpha \chi_{\{ X \geq X_{\mathrm{c}} \}} (X-X_{\mathrm{c}}),
\end{equation*}
with $v_0 = 1.76\times 10^{-3}\, \rm m/s$, $\breve{X} = 3.87\, \rm kg/m^3$, $\eta = 3.58$, $X_{\mathrm{c}} = 5\, \rm kg/m^3$ and $\alpha = 0.2\,\rm m^2/s^2$. 
Other parameters are $\rho_X = 1050\, \rm kg/m^3$, $\rho_L = 998\, \rm kg/m^3$, $g = 9.81 \,\rm m/s^2$, and $B = 3\,\rm m$.
The soluble feed concentrations in both examples are 
\begin{equation} \label{eq:Sf}
\bSf(t)\equiv (6.00\times 10^{-3}, 9.00\times 10^{-4},0)^{\rm T}\, \rm kg/m^3.
\end{equation} 

For visualization purposes, we do not plot zero  numerical  concentrations above the surface, but fill this region with grey colour.

\subsection{Example~1: An SBR cycle} \label{subsec:numex1}

\begin{table}[t]
\caption{Example 1: Time functions for an SBR cycle. `Model' refers to either PDE~\eqref{finalmod} or ODE~\eqref{eq:ODEmodel}.\label{table:SBRcycle}} 

\smallskip 

\centering
{\small
\begin{tabular}{lccccc} \toprule 
Stage  & Time period [h]& $\Qf(t) [\rm m^3/h]$ & $\Qu(t) [\rm m^3/h]$ & $\Qe(t) [\rm m^3/h]$ & Model\\
\midrule 
Fill & $0\leq t<1$ & 790 & 0& 0 & PDE\\
React & $1\leq t<3$ & 0 & 0& 0 & ODE\\
Settle & $3\leq t<5$ & 0 & 0& 0 & PDE\\
Draw & $5\leq t<5.5$ & 0 & 0& 1570 & PDE\\
Idle & $5.5\leq t<6$ & 0 & 10& 0 & PDE\\
\bottomrule 
\end{tabular}} 
\end{table}%

 \begin{figure}[t] 
\centering 
\begin{tabular}{cc}
\includegraphics[scale=0.44]{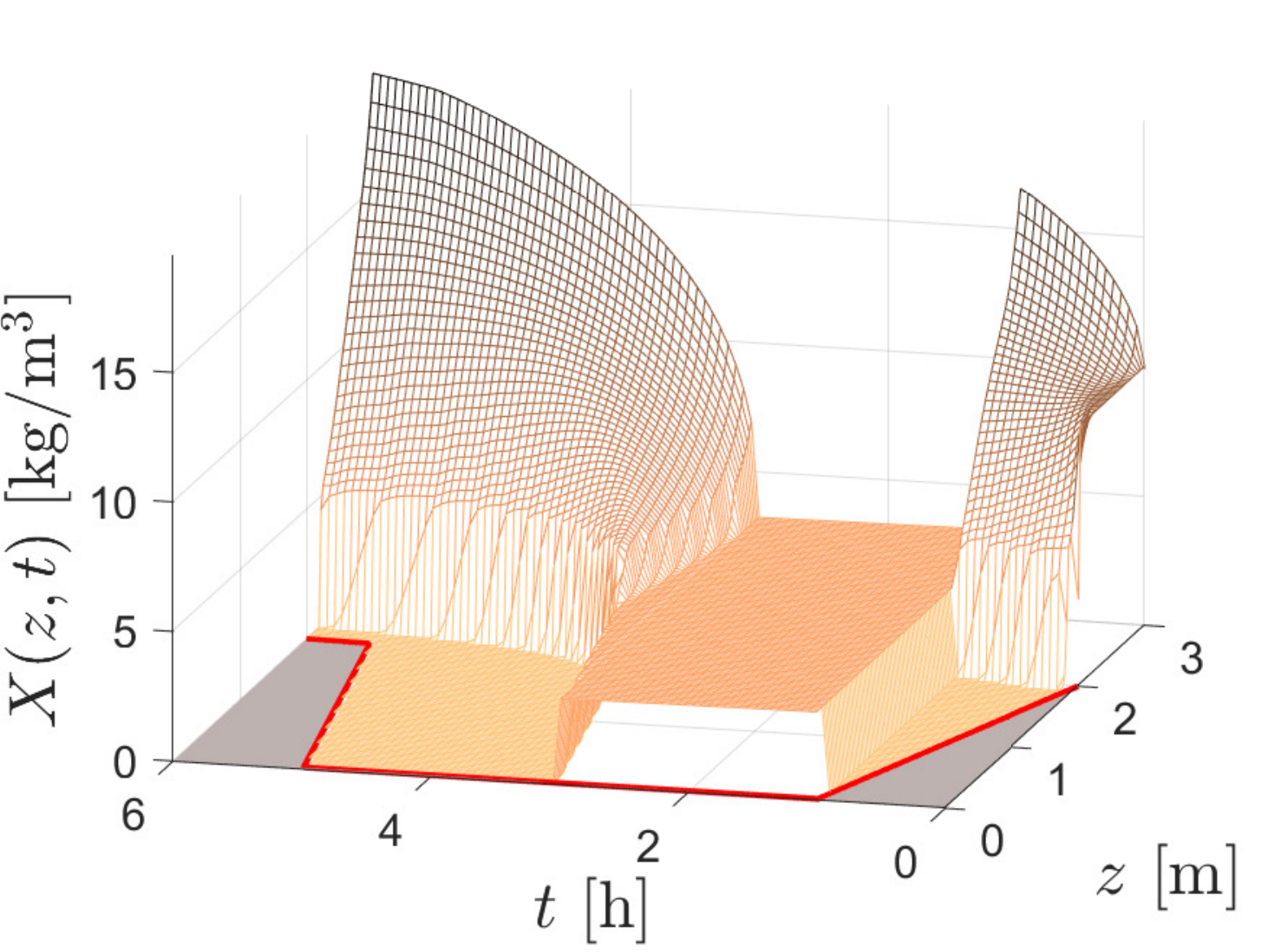} &
\includegraphics[scale=0.44]{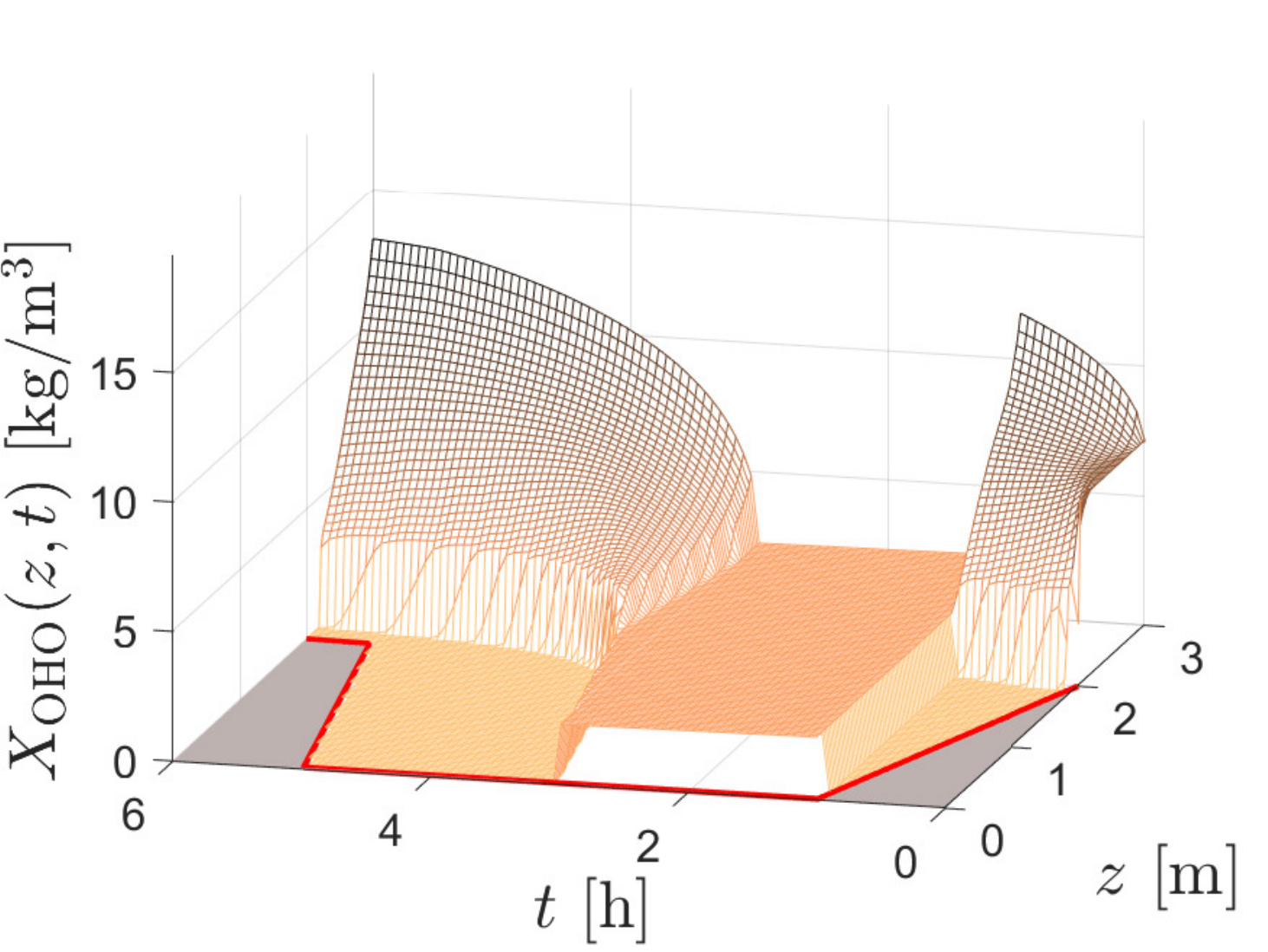} \\
\includegraphics[scale=0.44]{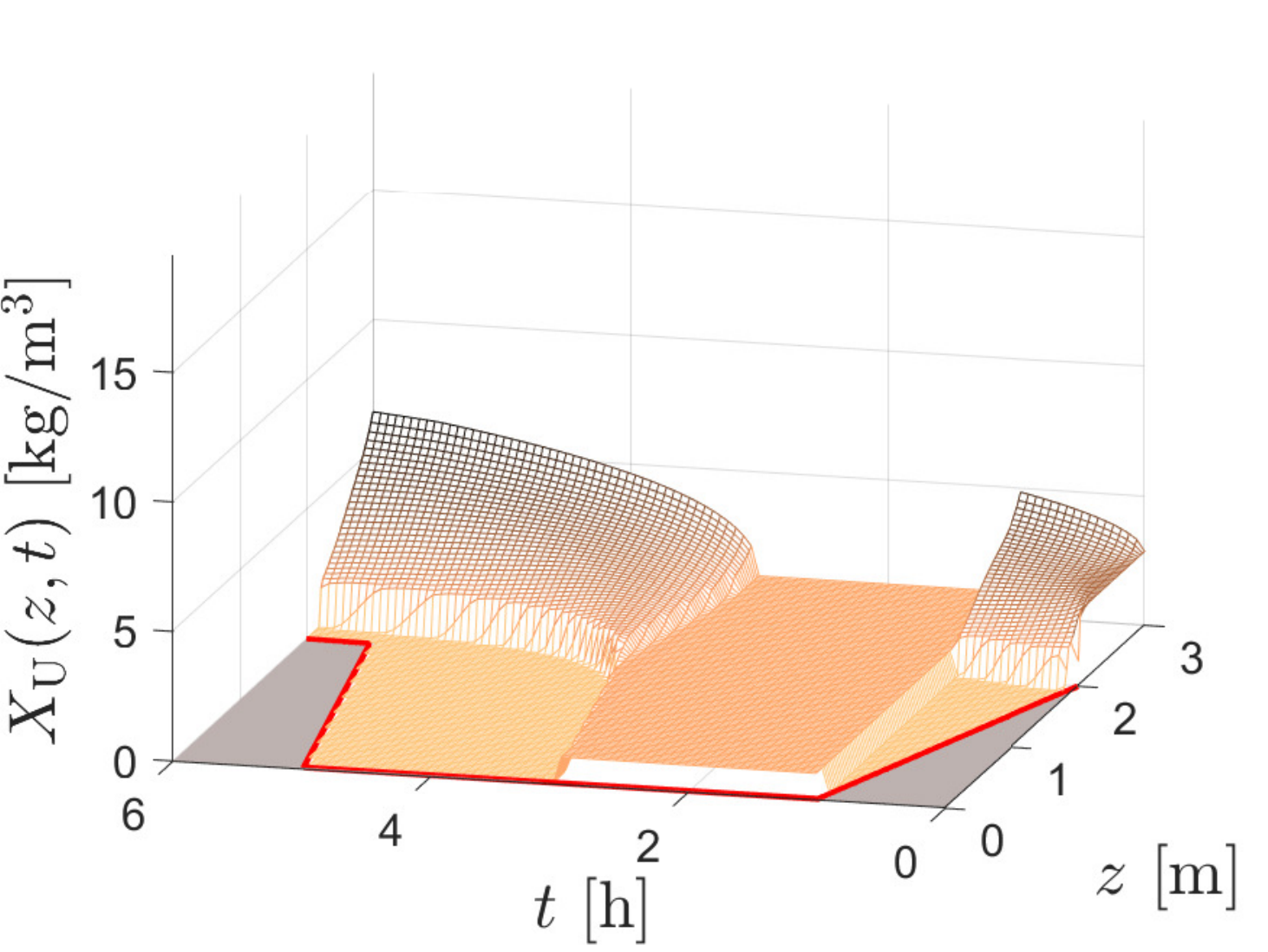} &
\includegraphics[scale=0.44]{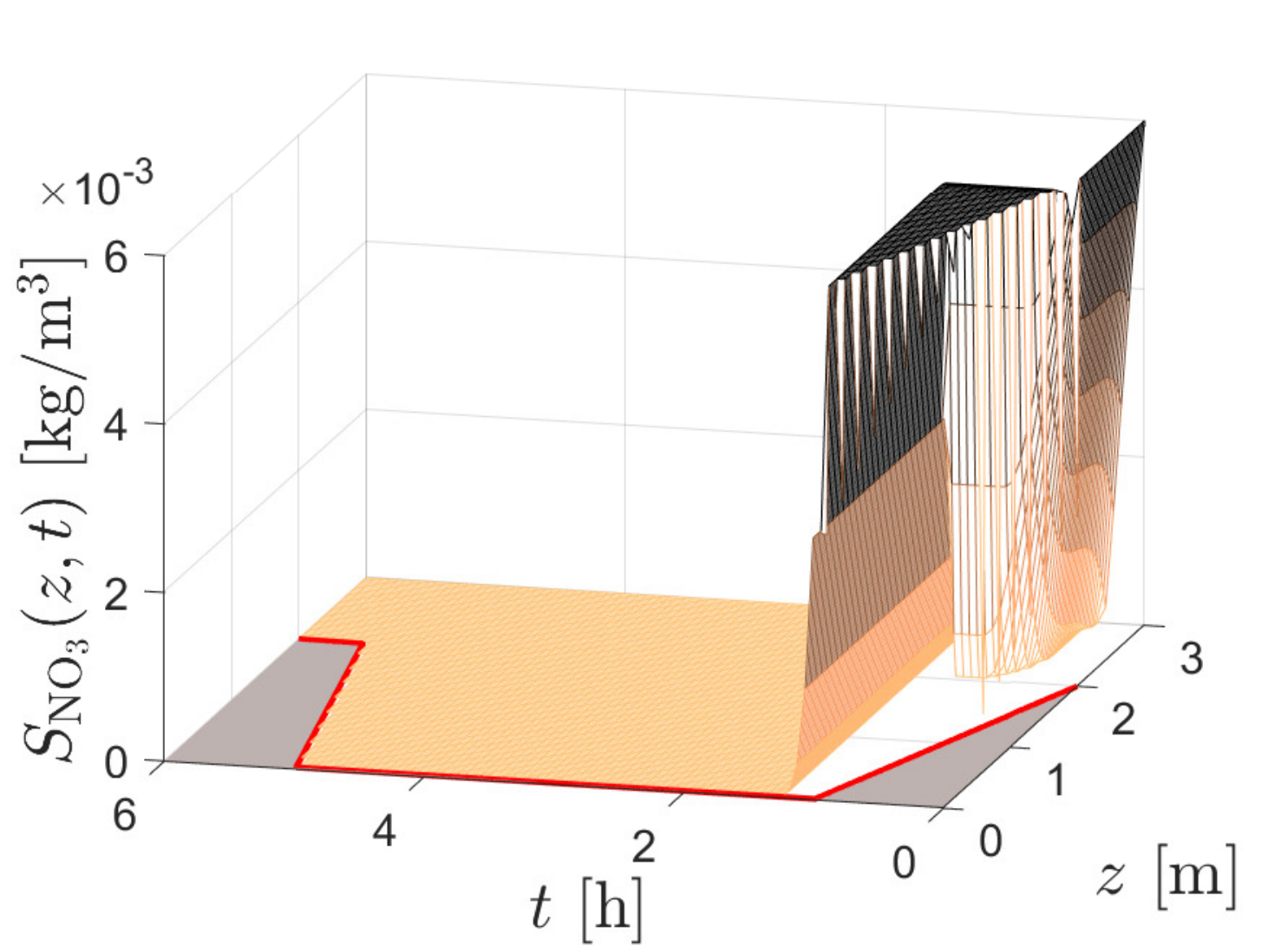}\\
\includegraphics[scale=0.44]{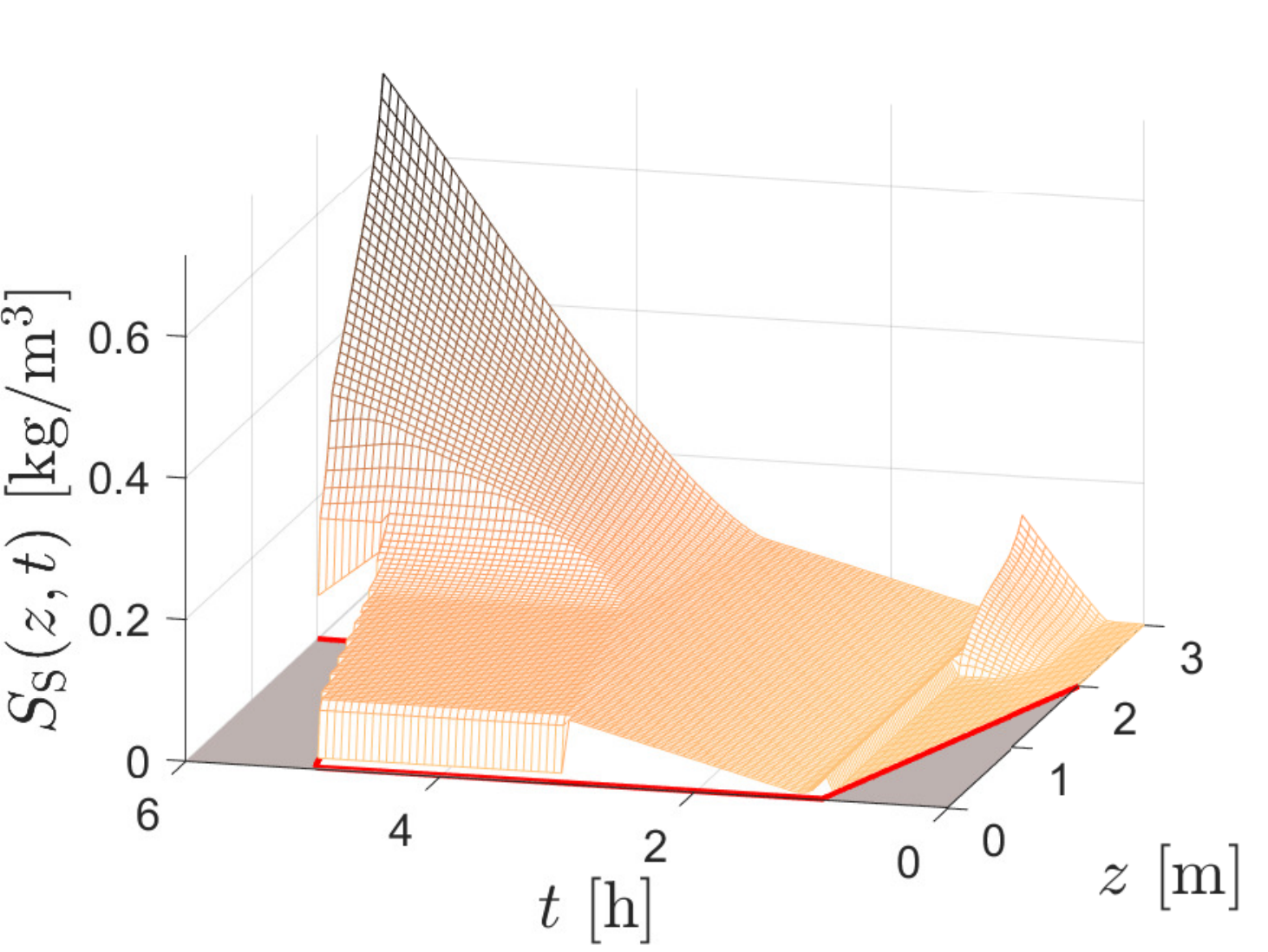} &
\includegraphics[scale=0.44]{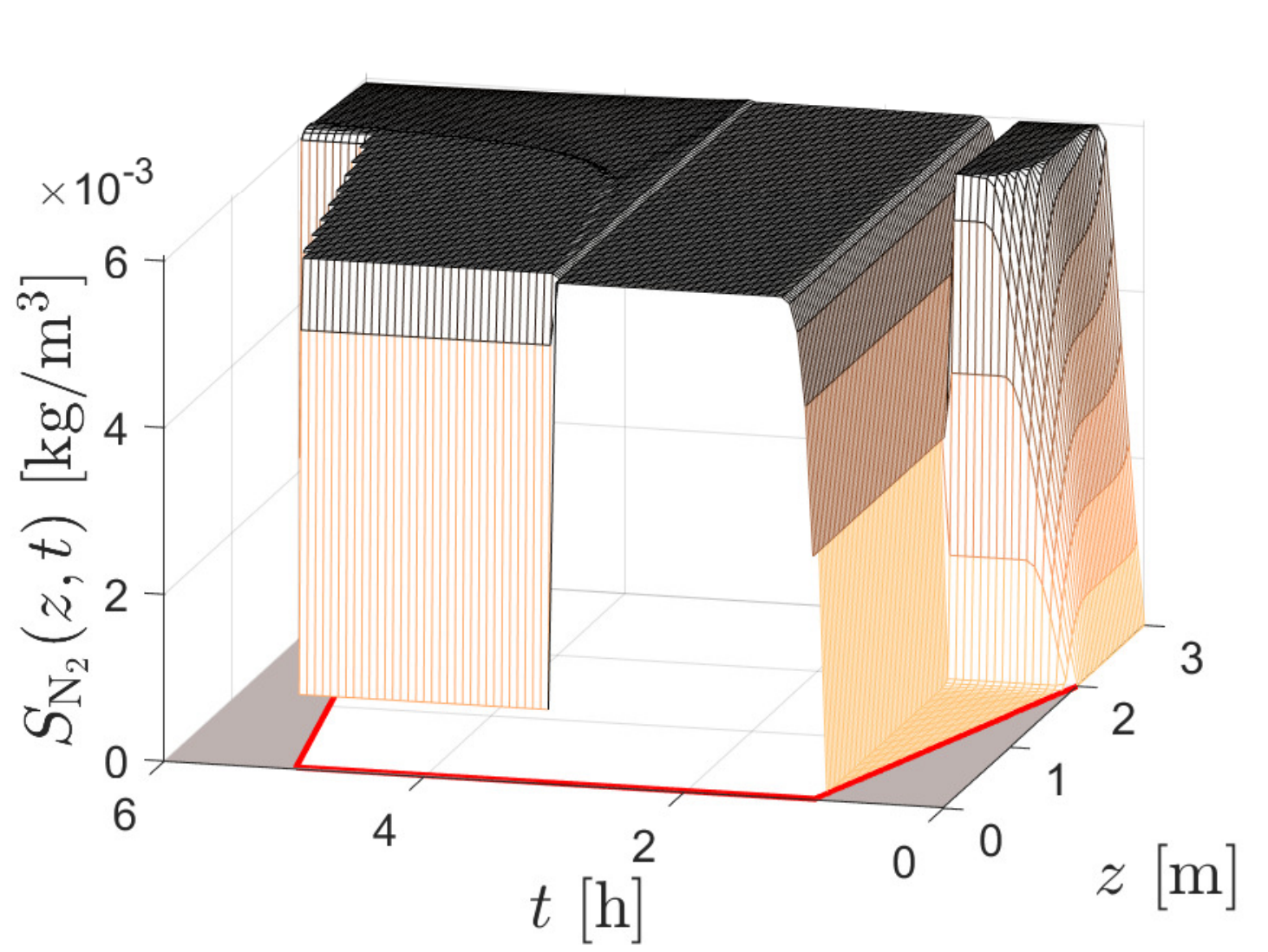}\\
\end{tabular}
 \caption{Example 1: Simulation results during $T = 6$ hours.
{The red lines show the location of the surface and the region above the surface is filled with grey colour.}} \label{fig:sim1}
\end{figure}%

A cylindrical tank with cross-sectional area $A = 400\,\rm m^2$ is simulated during~$6\,$h.
The lengths of the five stages are chosen primarily for illustration; see Table~\ref{table:SBRcycle}.
The initial concentrations are $\bC^0(z) = X^0(z) (5/7, 2/7)^{\mathrm{T}}$, where 
   \begin{alignat*}{3} 
      & X^0(z) =   0    \,{\rm kg/m^3},  \quad && \bS^0(z) =  \boldsymbol{0}   \,{\rm kg/m^3} 
      \quad  && \text{if  $z<2.0 \,{\rm m}$,} \\
      & X^0(z) =  10\,{\rm kg/m^3}   \quad && \bS^0(z) 
      = (6\times 10^{-3}, 9\times  10^{-4}, 0)^{\rm T}\,{\rm kg/m^3} 
      \quad &&  \text{if $z\geq 2.0\,\rm m $.} 
    \end{alignat*}
No biomass is fed to the tank; $\bCf(t)\equiv\bzero$.
Figure~\ref{fig:sim1} shows the simulation results. 
The reactions converting $\rm NO_3$ to $\rm N_2$ start immediately and are fast. (The downwards-pointing peaks in the $S_{\mathrm{NO}_{\mathrm{3}}}$ plot arise since we do not plot  zero concentrations above the surface.)
A short time after the react stage has started at $t=1\,$h, all $\rm NO_3$ has been consumed.
During this short time period, $S_{\rm S}$ decreases slightly when there is still sufficient $\rm NO_3$, but then increases during the react stage because of the decay of heterotrophic organisms.

\subsection{Example 2} \label{subsec:numex2}

\begin{figure}[t] 
\centering 
\includegraphics[scale=1.25]{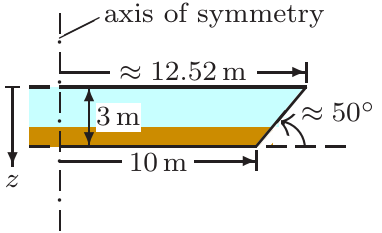} 
\caption{Example 2:  Schematic of the truncated cone. \label{fig:trunccone}} 
\end{figure}%

\begin{table}[t]
\caption{Example 2. Schematic of the truncated cone and time functions for the simulation. `Model' refers to either PDE~\eqref{finalmod} or ODE~\eqref{eq:ODEmodel}.\label{table:Ex2}}

\smallskip 

\centering
{\small  \begin{tabular}{cccccc} \toprule
Time period $[\mathrm{h}]$ & $\Xf(t) [\rm kg/m^3]$& $\Qf(t) [\rm m^3/h]$& $\Qu(t) [\rm m^3/h]$& $\Qe(t) [\rm m^3/h]$ & {Model} \\
\midrule
$0\leq t<1$ & 5 & 790 & 0& 0 & PDE\\
$1\leq t<2$ &0  &  0 & 100& 0 & ODE\\
$2\leq t<3$ &0  &  0 & 0& 100 & ODE\\ 
$3\leq t<5$ & 5  &  100 & 0& 0 & PDE\\
$5\leq t<6$ &0  &  0 & 0& 790 & PDE\\
\bottomrule
\end{tabular}}  
\end{table}%

We now choose a truncated cone  (cf.\ Figure~\ref{fig:trunccone}) of  the same volume $1200\,\rm m^3$ as the cylinder of  Example~1 and demonstrate what the numerical scheme can handle during extreme cases of fill and draw when solids concentrations are positive at the surface.
We use the same initial data as in Example~1 but with  $\bar{z}(0)\approx 1.8429\,$m to obtain the same initial volume of mixture as in Example~1.
The fill and draw periods are specified in Table~\ref{table:Ex2}.
The feed concentrations of the substrates are given by~\eqref{eq:Sf} and those of the biomass by
 $\bC_{\rm f}(t) = X_{\rm f}(t) (5/7, 2/7)^{\mathrm{T}}$, where the piecewise constant 
  function $X_{\rm f}(t) $ follows from  Table~\ref{table:Ex2}. 

 \begin{figure}[t] 
\centering 
 \begin{tabular}{cc}
\includegraphics[scale=0.44]{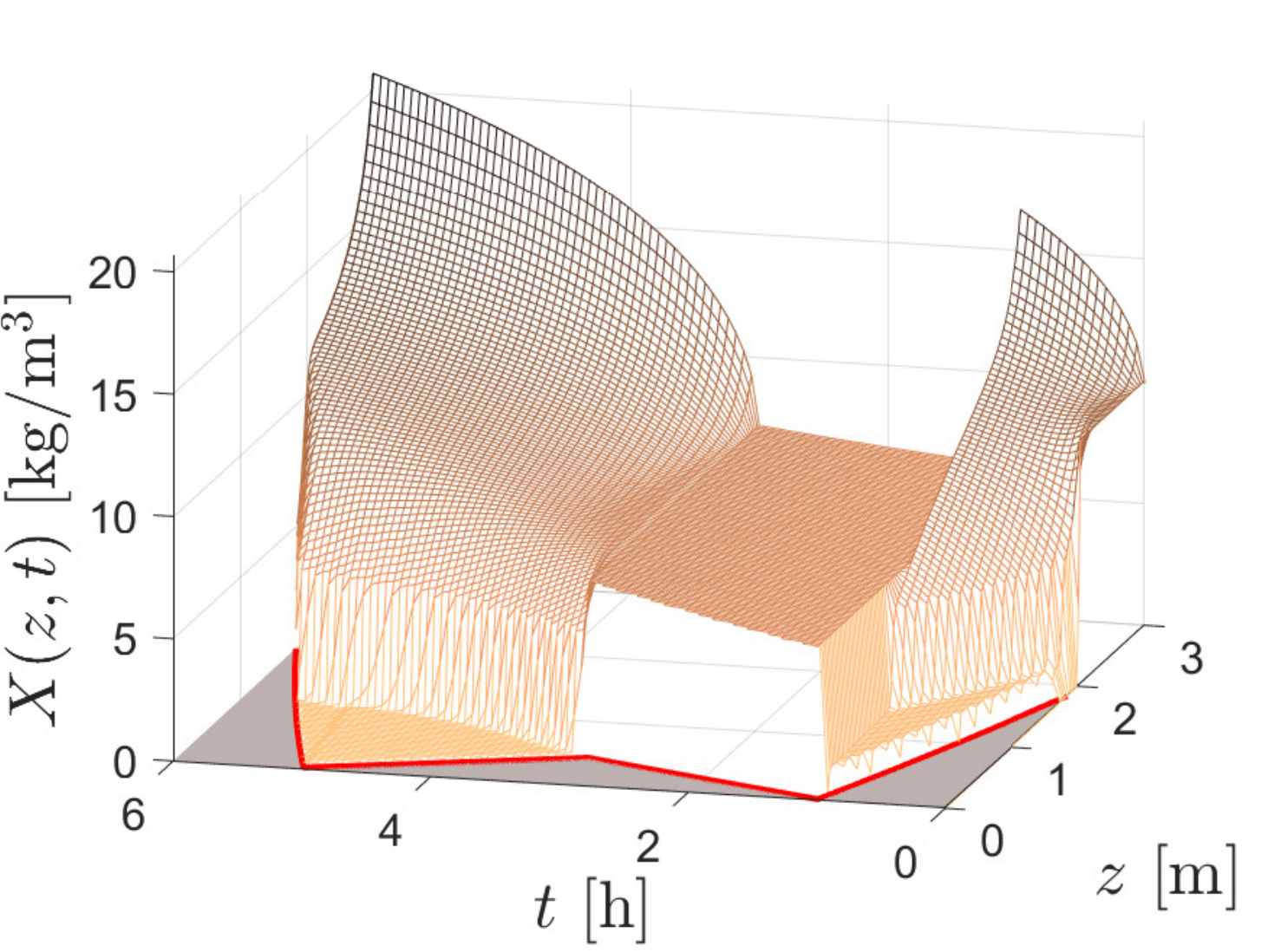} &
\includegraphics[scale=0.44]{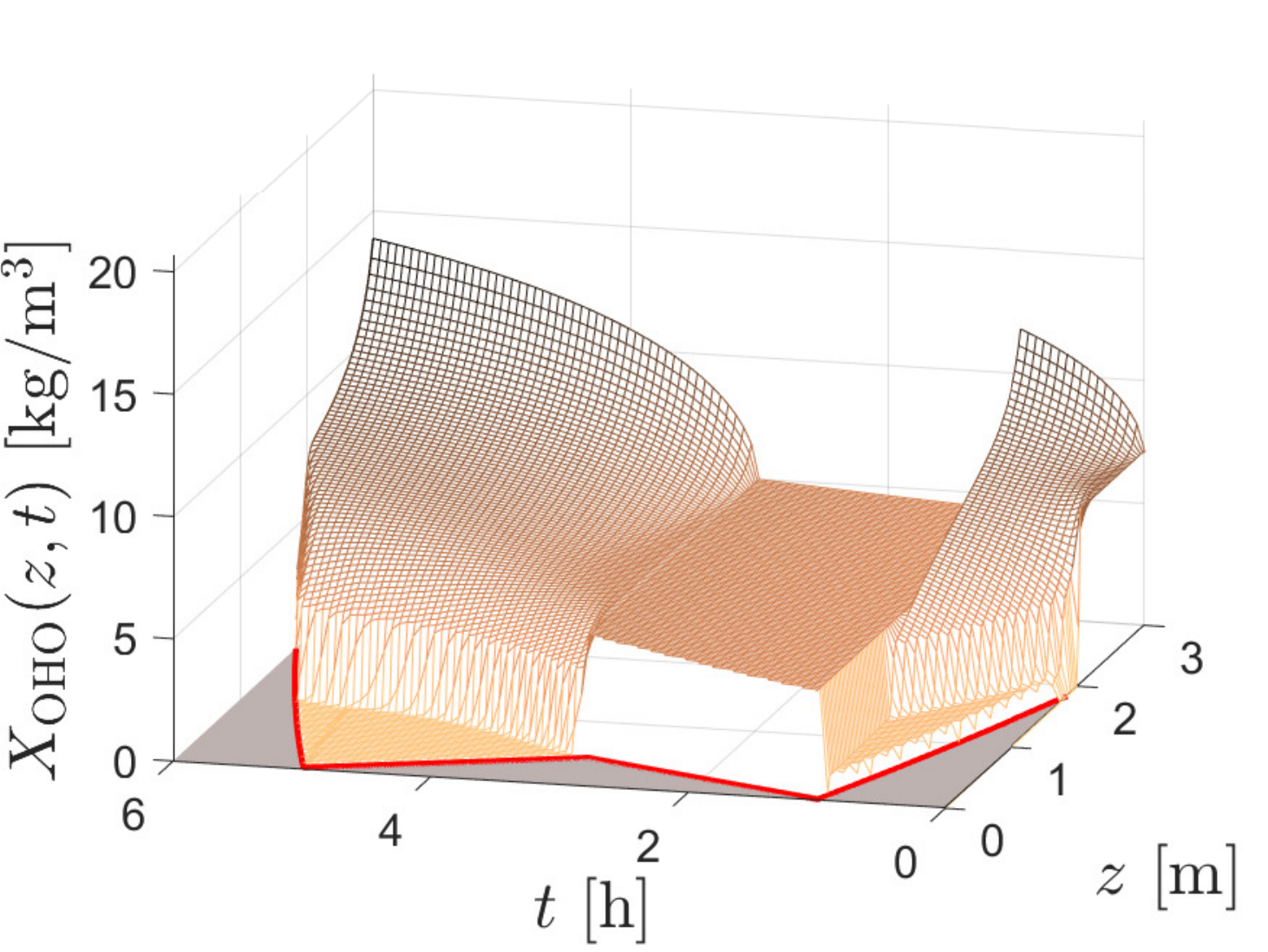} \\
\includegraphics[scale=0.44]{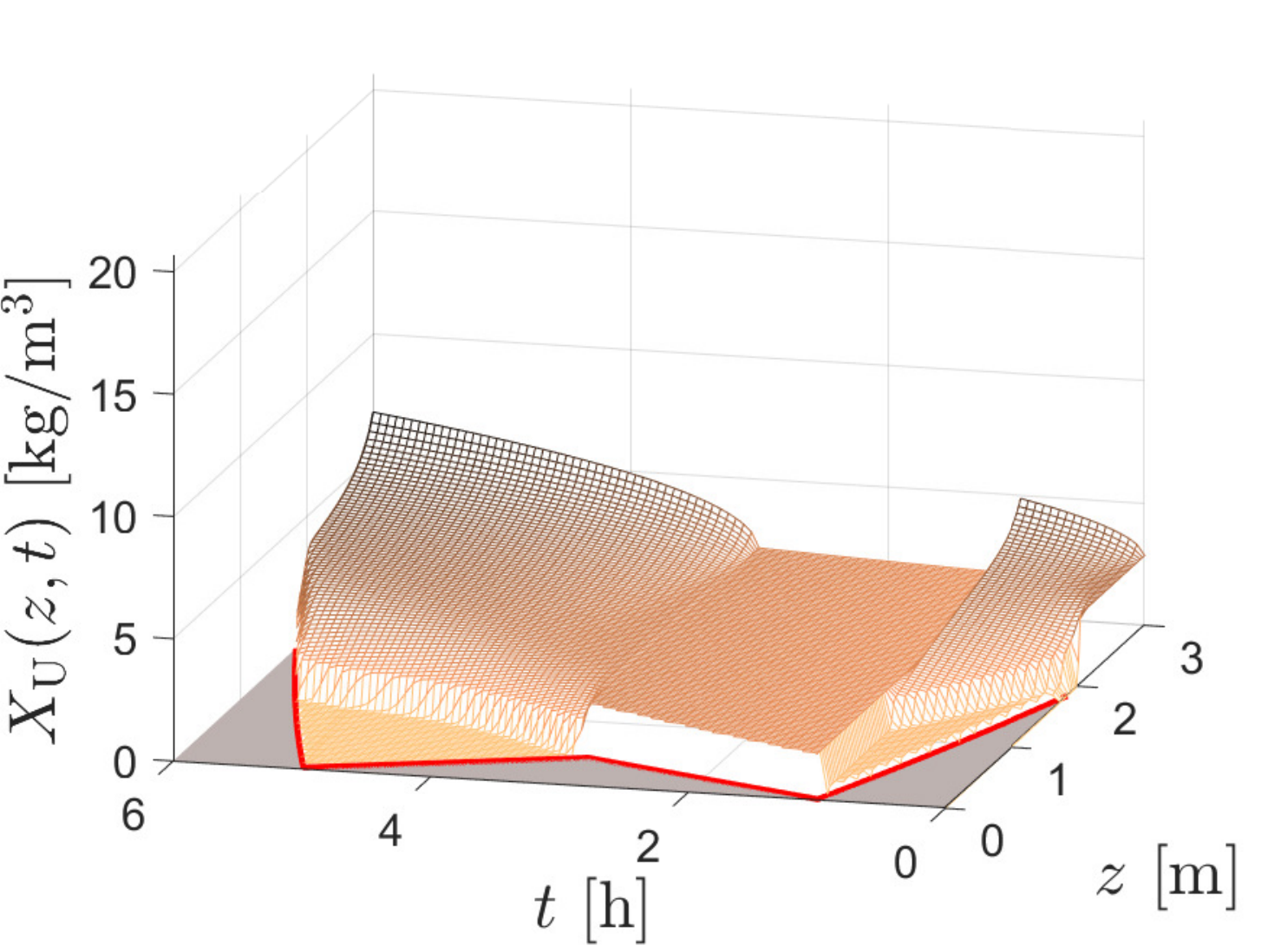} &
\includegraphics[scale=0.44]{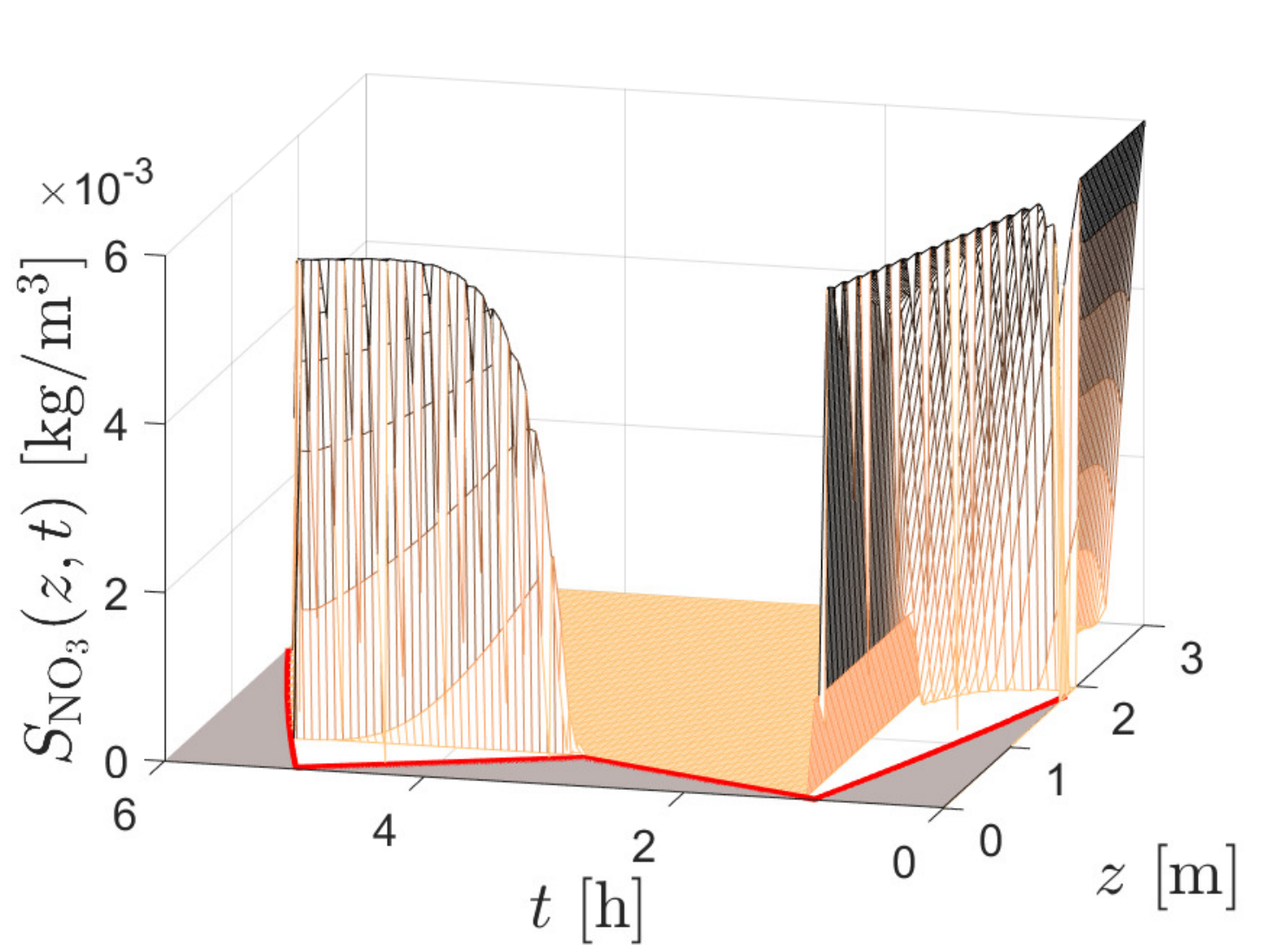}\\
\includegraphics[scale=0.44]{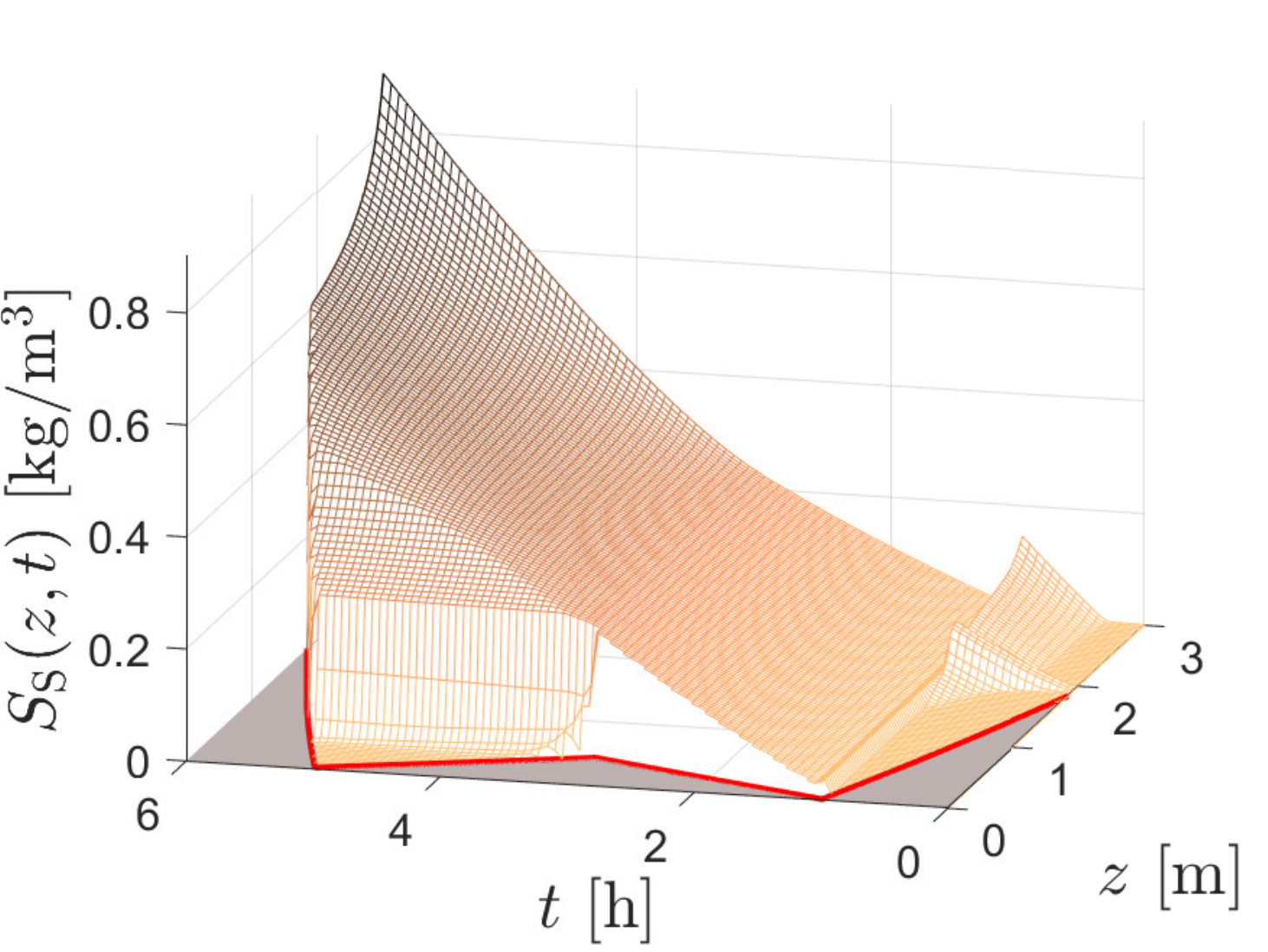} &
\includegraphics[scale=0.44]{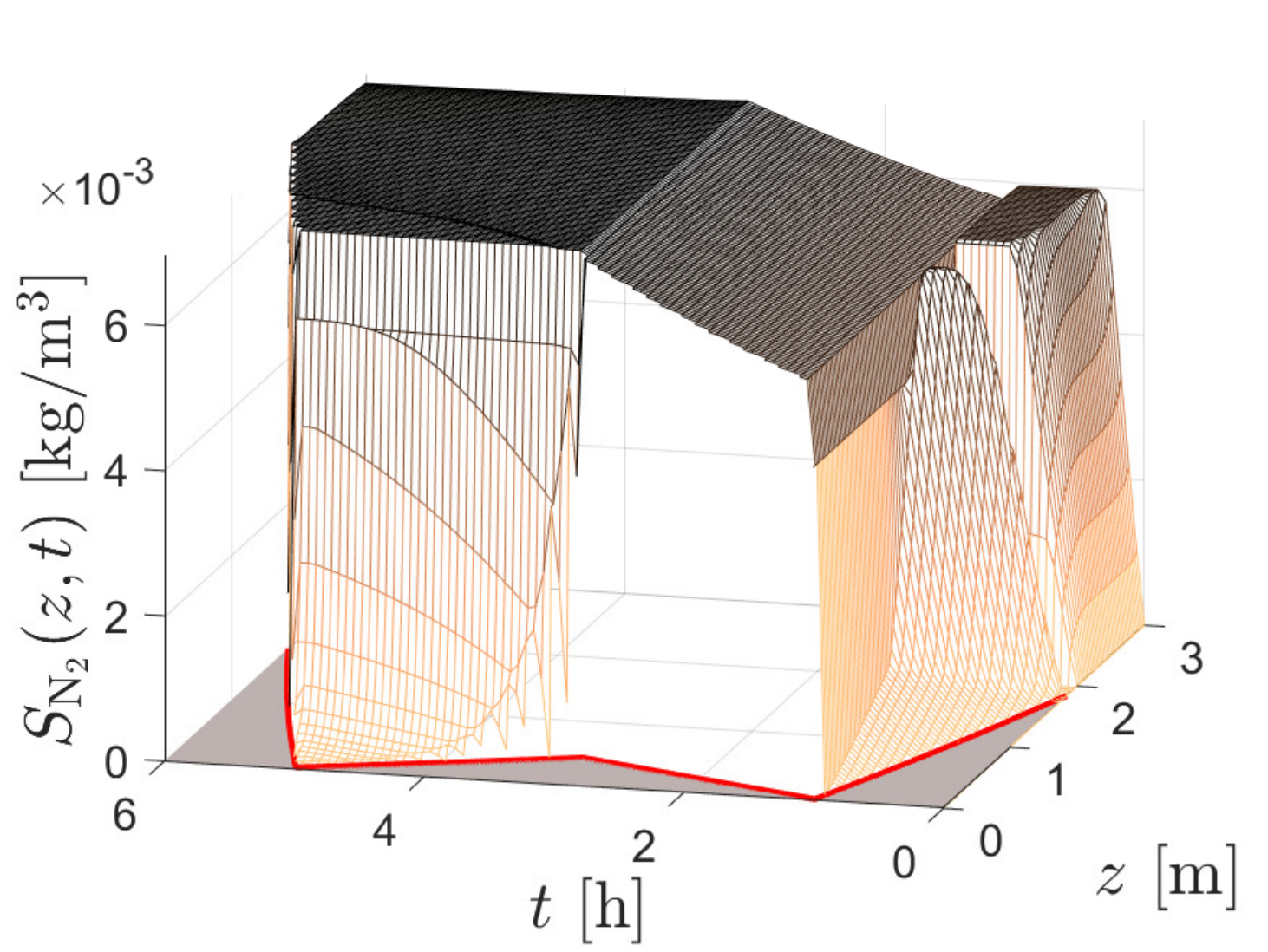}\\
 \end{tabular}
\caption{Example 2: Simulated results during $T = 6$ hours.}  \label{fig:sim2}
\end{figure}%

\begin{figure}[t] 
\centering 
\begin{tabular}{cc}
\includegraphics[scale=0.40]{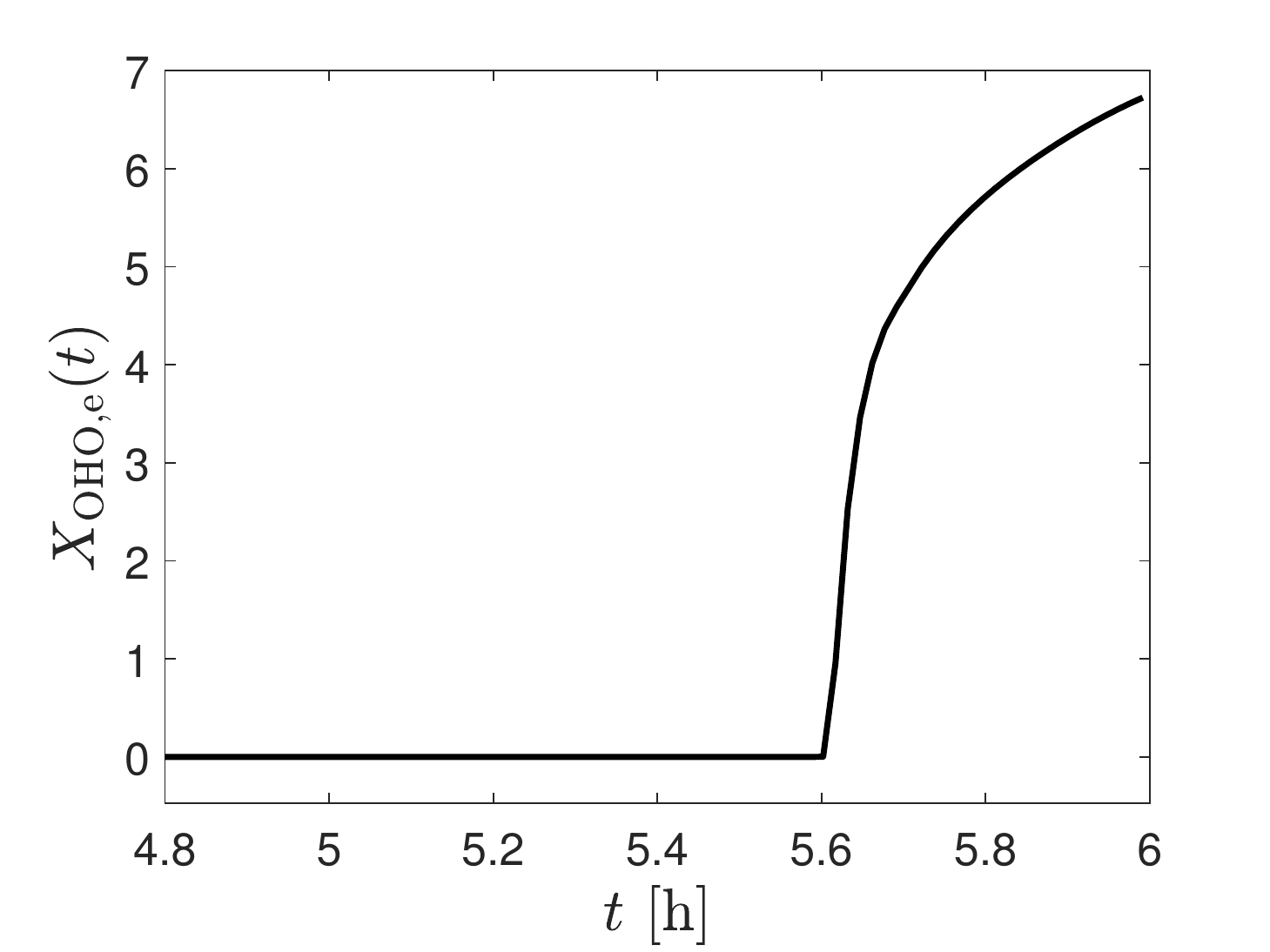} &
 \hspace{-0.5cm} \includegraphics[scale=0.40]{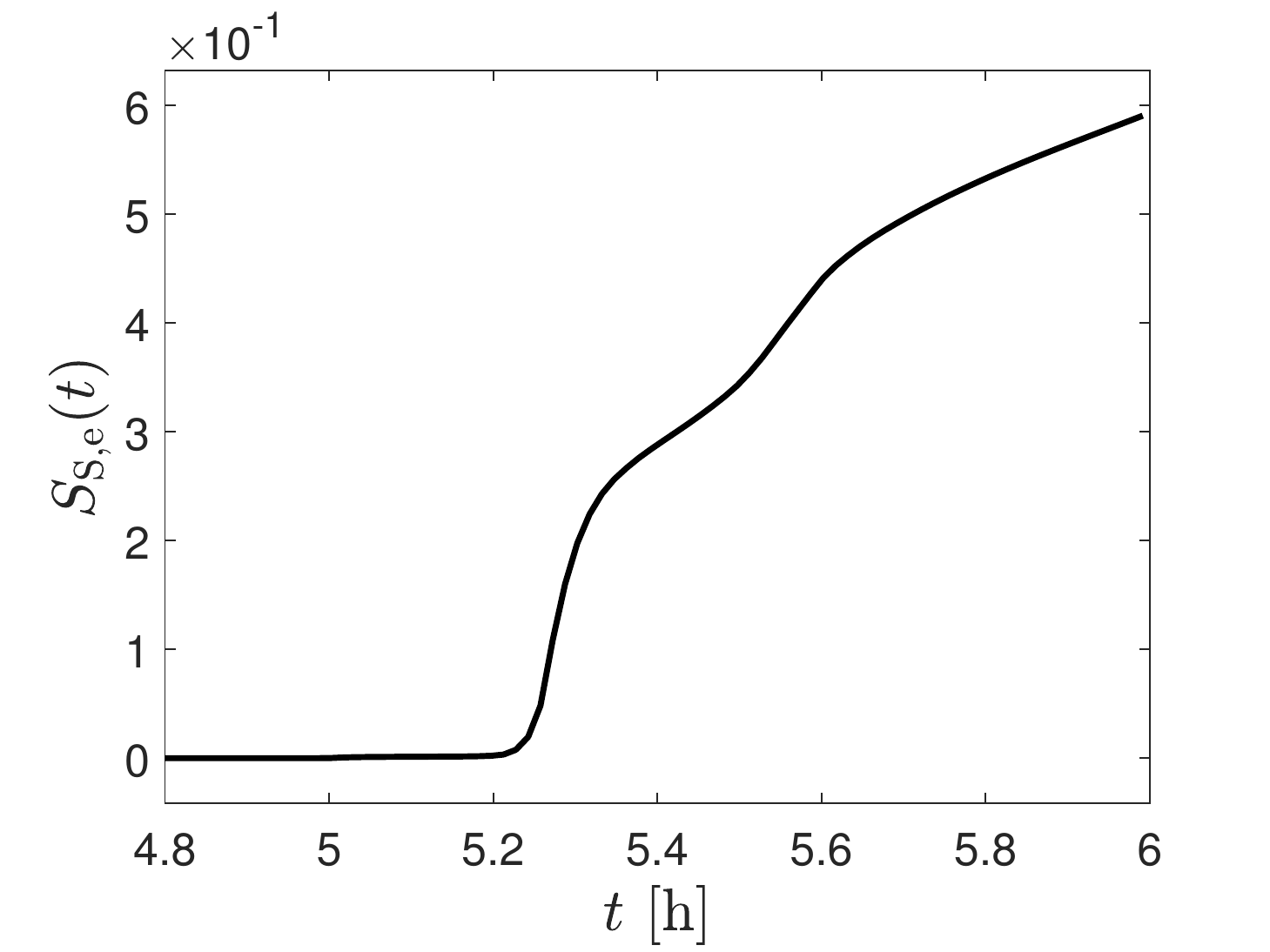} \\
\includegraphics[scale=0.40]{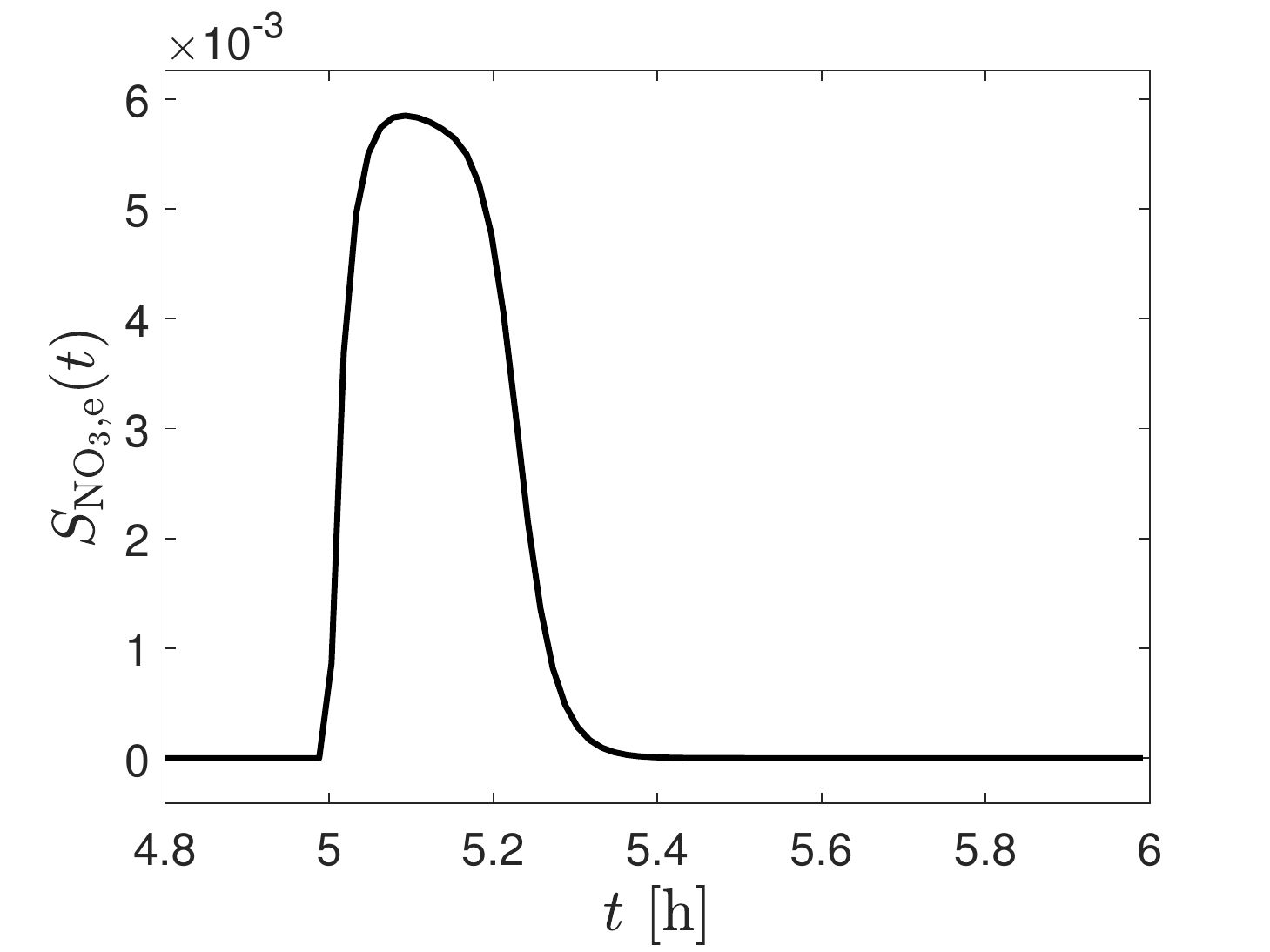} &
 \hspace{-0.5cm} \includegraphics[scale=0.40]{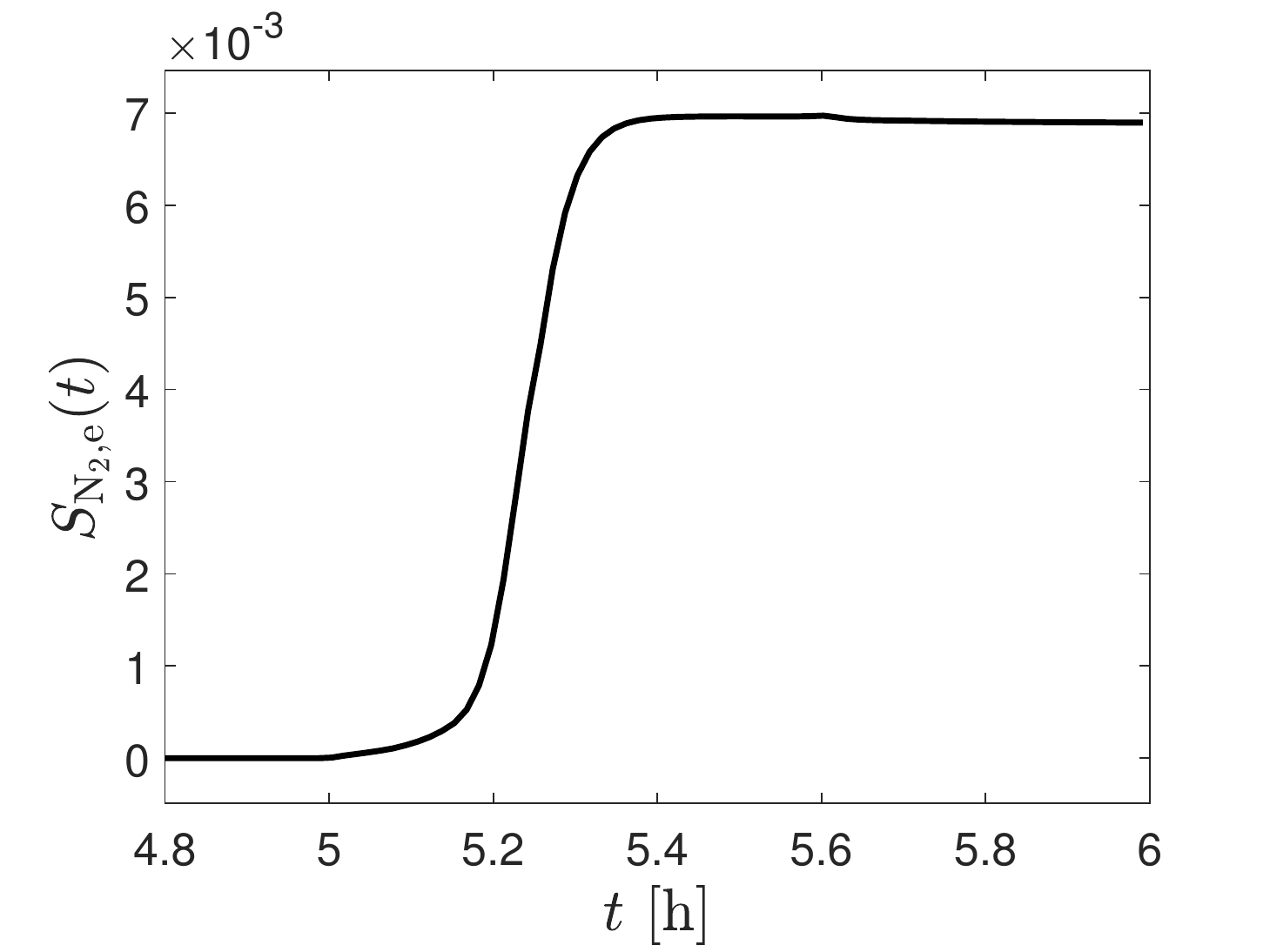}\\
\end{tabular}
\caption{Example 2: Simulated effluent concentrations, all in $\mathrm{kg}/\mathrm{m}^3$,  during 
  $t \in [4.8 \, \mathrm{h}, 6 \, \mathrm{h}]$ obtained by the discretization where the number of computational cells within the tank is $100$.}\label{fig:effluent}
\end{figure}%
  
Figure~\ref{fig:sim2} shows the simulated concentrations.
During the first hour, there is a discontinuity in the solids concentration $X$ rising with a lower speed than the surface.
Then full mixing occurs during two hours and the surface is lowered because of the outlet flows at the bottom and top.
At $t=3\,\rm h$, the mixing stops and the solids settle again.
During $3\,{\rm h}<t<5\,{\rm h}$, the tank is filled up again with solids and substrates.
The solids feed concentration $\Xf=5\,\rm kg/m^3$ is the same as during the first hour, but now the feed flow $\Qf$ is much lower, and hence the mass flow much lower.
The result is a very low concentration~$X$ below the surface during $3\,{\rm h}<t<5\,{\rm h}$.
Since also $X_{\rm OHO}$ is low, there are  hardly any reactions and most of the fed $\rm NO_3$ remains in the mixture above the sludge blanket of the solids.
At the surface level around $t=3\,$h, there is also biomass present and a high production of $\rm N_2$ occurs.
However, the sludge blanket drops and the high concentration of $\rm N_2$ remains at this height until it is extracted through the effluent during the last hour.
The latter is shown in Figure~\ref{fig:effluent}, which also shows that solids are extracted.

\section{Conclusions} \label{sec:concl}

A general model of multi-component reactive settling of flocculated particles given by a quasi-one-dimensional PDE system with moving boundary~\eqref{eq:model}, \eqref{eq:modelterms} is introduced.
Fill and draw of mixture at the moving surface can be made at any time and a specific application is the SBR process.
The unknowns are concentrations of solids and soluble substrates, and the PDE model can (via its reaction terms) be combined with well-established models for the biochemical reactions in wastewater treatment.

The moving boundary can be precomputed with the ODE~\eqref{eq:zbarprime} containing the volumetric flows in and out of the tank.
{The positivity-preserving numerical scheme of~\cite{bcdp_part2} is used to simulate the SBR process when denitrification occurs in a tank with either constant or a varying cross-sectional area and when extreme cases of fill and draw occur. Indication of the convergence of the numerical scheme as the mesh size is reduced are provided in~\cite{bcdp_part2}.}

With the present model, investigations and optimization of the SBR process can be made, and the usage of several SBRs coupled in series or in parallel with synchronized stages so that, for instance, a continuous stream of effluent of certain quality is obtained.
Furthermore, more accurate comparisons are possible between SBRs and continuously operated SSTs, since one or the other may be preferred depending on the plant size and other practical considerations~\cite{droste}.

\section*{Acknowledgements}
RB~is  supported by ANID (Chile) through 
 projects Centro de Modelamiento Matem\'{a}tico  (BASAL pro\-jects ACE210010 and FB210005); 
  Anillo project ANID/PIA/ACT210030;  CRHIAM, project   ANID/FONDAP/15130015;  and  Fon\-de\-cyt project 1210610. SD~acknowledges support from the Swedish Research Council (Vetenskapsr\aa det, 2019-04601). RP  is supported by ANID scholarship ANID-PCHA/Doctorado Nacional/2020-21200939.

\bibliography{ref_copy}

\end{document}